\magnification 1200
\voffset 20pt

\

\font\bigbold=cmbx12 scaled 1000 
\font\smallit=cmti8

\font\Bbb=msbm10
\def\BBB#1{\hbox{\Bbb#1}}
\def\C{{\BBB C}}
\def\Z{{\BBB Z}}
\def\Q{{\BBB Q}}

\font\Bbbs=msbm10 scaled 800
\def\BBBS#1{\hbox{\Bbbs#1}}
\def\CS{{\BBBS C}}
\def\ZS{{\BBBS Z}}

\def\o{{\bf 1}}

\def\z{\left[ z_1^{-1} \delta \left( {z_2 \over z_1} \right) \right]}
\def\zd{\left[ z_1^{-1}{\d \over \d z_2}\delta \left( {z_2 \over z_1} \right) \right]}

\def\CC{{\cal C}}
\def\D{{\cal D}}
\def\H{{\cal H}}
\def\L{{\cal L}}
\def\S{{\cal S}}
\def\R{{\cal R}}
\def\U{{\cal U}}
\def\F{{\cal F}}
\def\hR{{\widehat \R}}
\def\hH{{\widehat \H}}
\def\hT{{\widehat T}}

\def\bV{{\overline V}}
\def\tV{{\widetilde V}}
\def\tW{{\widetilde W}}
\def\d{\partial}
\def\ot{\otimes}
\def\td{{\tilde d}}
\def\lcol{\hbox{{\kern+0.02em}{\bf :}\kern-.02em}}
\def\rcol{\hbox{\kern-.02em{}{\bf :}}{\kern+0.02em}}

\def\Der{{\rm Der \;}}
\def\End{{\rm End}}
\def\mod{{\rm mod \;}}
\def\sll{{sl}}
\def\gll{{gl}}
\def\Ind{{\rm Ind}}
\def\Span{{\rm Span}}
\def\id{{\rm id}}
\def\tr{{\rm tr}}

\def\z{\left[ z_1^{-1} \delta \left( {z_2 \over z_1} \right) \right]}

\def\zd{\left[ z_1^{-1}  {\d \over \d z_2}
\delta \left( {z_2 \over z_1} \right) \right]}

\def\isot{0.1}
\def\odz{0.2}
\def\odo{0.3}
\def\ode{0.4}
\def\sole{0.5}
\def\oap{0.6}
\def\ovo{0.7}
\def\opi{0.8}
\def\opii{0.9}

\def\Lieb{1.1}
\def\charf{1.2}

\def\comm{2.1}
\def\vla{2.2}
\def\Y{2.3}
\def\done{2.4}
\def\dzer{2.5}
\def\doomr{2.6}
\def\dozmr{2.7}
\def\dzzmr{2.8}

\def\doo{3.1}
\def\dzm{3.2}
\def\dom{3.3}
\def\drm{3.4}
\def\rais{3.5}
\def\Ydzm{3.6}
\def\Ydom{3.7}
\def\Ydrm{3.8}
\def\relm{3.9}
\def\drkm{3.10}
\def\dooz{3.11}
\def\dooo{3.12}
\def\dozz{3.13}
\def\dozo{3.14}
\def\dozt{3.15}
\def\dzzz{3.16}
\def\dzzo{3.17}
\def\dzzt{3.18}
\def\fdrm{3.19}
\def\frkm{3.20}
\def\fmm{3.21}
\def\Wact{4.1}
\def\spec{4.2}

\def\vdom{5.1}
\def\nopm{5.2}
\def\vom{5.3}
\def\otw{5.4}
\def\opw{5.5}
\def\szp{5.6}

\def\diffeq{6.1}
\def\eqw{6.2}
\def\eqA{6.3}
\def\zint{6.4}
\def\cons{6.5}
\def\Pnz{6.6}
\def\Poz{6.7}
\def\wj{6.8}

\def\PPsn{A.1}
\def\PTsn{A.2}
\def\FGrel{A.3}

\def\voa{2.1}
\def\inv{2.2}
\def\vl{2.3}

\

\centerline{\bigbold Differential equations in vertex algebras and 
simple modules}
\centerline{\bigbold for the Lie algebra of vector fields on a torus.}

\

\centerline{\bf Yuly Billig${}^1$, Alexander Molev${}^2$ and Ruibin Zhang${}^2$}

{
\smallit
\narrower
\noindent
${}^1$Carleton University, School of Mathematics and Statistics,
1125 Colonel By Drive, Ottawa, \hfill
\break
K1S 5B6, Canada \hfill
\break
\noindent
${}^2$University of Sydney, School of Mathematics and Statistics,
Sydney, NSW 2006, Australia
}

\

\

{\bf Abstract.} 
We study irreducible representations for the Lie algebra of vector fields
on a 2-dimensional torus constructed using the generalized Verma modules.
We show that for a certain choice of parameters these representations remain 
irreducible when restricted to a loop subalgebra in the Lie algebra of vector 
fields.
We prove this result by studying vertex algebras associated with the Lie algebra of vector fields on a torus and solving non-commutative differential equations that we derive using the vertex algebra technique.

\

{\bf Mathematics Subject Classification (2000):} 17B66, 17B69.

{\bf Key words:} Lie algebra of vector fields, vertex algebras, generalized Verma modules.

\

\

{\bf 0. Introduction.}

\

In recent years a new area in representation theory has emerged -- 
the theory of bounded modules for infinite-dimensional Lie algebras 
with a dense $\Z^n$-grading. The classical case $n=1$ includes 
Kac-Moody algebras and the Virasoro algebra, and the main tool in 
the representation theory in this case is the concept of a Verma module. For 
$n > 1$, however, there is no natural way to split 
$\Z^n \backslash \{ 0 \}$
into a positive and a negative parts, and the Verma module technique 
does not yield interesting irreducible modules.

 Nonetheless Berman and Billig [BB] showed that for a class of Lie 
algebras with a polynomial multiplication, the generalized Verma 
modules obtained by cutting $\Z^n$ with a hyperplane that intersects 
$\Z^n$ in a lattice of rank $n-1$, have irreducible quotients with 
finite-dimensional weight spaces. Billig and Zhao [BZ] extended this 
result to Lie algebras and modules with exp-polynomial multiplication 
and gave many examples illustrating this theorem.

 In the present paper we investigate representations of one of the most 
natural Lie algebras with a dense $\Z^2$-grading -- the Lie algebra of
vector fields on a 2-dimensional torus:
$$\D = \Der \C[t_0^{\pm 1}, t_1^{\pm 1}] .$$
This Lie algebra is $\Z^2$-graded by the eigenvalues of the adjoint action 
of
$$ d_0 = t_0 {\d \over \d t_0 } \quad {\rm and} \quad
d_1 = t_1 {\d \over \d t_1 } .$$
This grading is dense, i.e., its support is all of $\Z^2$. 

One class of representations of this Lie algebra with finite-dimensional 
weight spaces, the tensor modules, is well understood. These modules 
have a geometric realization as tensor fields on a torus.

 Here we will study another family of simple modules for $\D$ -- a class 
of bounded modules. Bounded modules are constructed using the general 
scheme given in [BB]. We consider a decomposition 
$$\D = \D_- \oplus \D_0 \oplus \D_+ $$
given by the eigenvalues of $d_0$.

 The degree zero part, $\D_0$, is isomorphic to a semidirect product of the Lie algebra
of vector fields on a circle with the abelian Lie algebra of functions
on a circle. Tensor modules for $\D_0$ are parametrized by triples
$\alpha, \beta, \gamma \in \C$. The module $T = T(\alpha, \beta,\gamma)$
has a basis $\{ v(n) | n \in \gamma + \Z \}$, and the action of $\D_0$
is given by the formulas
$$\eqalign{
(t_1^m d_1) v(n) &= (n - \alpha m) v(n+m) , \cr
(t_1^m d_0) v(n) &= \beta v(n+m) . \cr}$$

We use these $\D_0$-modules to construct the generalized Verma modules
$M(\alpha, \beta, \gamma)$ for the Lie algebra $\D$. Although      
$M(\alpha, \beta, \gamma)$ has infinite-dimensional weight spaces below
the top $T$, the general result of [BB] implies that its irreducible
quotient $L(\alpha, \beta, \gamma)$ has finite-dimensional weight spaces.

It is natural to ask for the dimensions of the weight spaces 
in $L(\alpha, \beta, \gamma)$. Here we give an explicit realization and
a character formula for these modules when $\alpha \not\in\Q$ and 
$\beta = -{\alpha (\alpha+1) \over 2}$. We show that in this case the picture
is analogous to the situation with the basic module for affine Kac-Moody 
algebras. 

Recall that an affine Kac-Moody algebra has an infinite-dimensional (principal)
\break
Heisenberg subalgebra, and its basic module remains irreducible when restricted
to the Heisenberg subalgebra [LW], [K1]. This yields a realization of the basic
module as a Fock space, and the action of the affine algebra is given by means
of the vertex operators that are constructed as solutions of certain 
operator-valued differential equations.

For the class of modules for the Lie algebra of vector fields that we consider
here, the role of the Heisenberg subalgebra is played by the loop subalgebra
$\L = \C [t_0, t_0^{-1} ] \otimes \sll_2 (\C) \subset \D$.

Our main result, Theorem 1.2, states that for $\alpha \not\in \Q$, 
$\beta = -{\alpha (\alpha+1) \over 2}$, the module $L(\alpha, \beta, \gamma)$
remains irreducible when restricted to the loop subalgebra $\L$.
We shall show that $L(\alpha, \beta, \gamma)$ can be realized
as an induced module for $\L$:
$$L(\alpha, \beta, \gamma) \cong U(\L_-) \otimes T, \eqno{(\isot)}$$
from which we immediately get a formula for the dimensions of the weight spaces.

To establish this result, we develop new methods that seem to be 
interesting in  themselves. 
The main tool we use is the machinery of vertex algebras. We consider
for the Lie algebra $\D$ the universal enveloping vertex algebra $V_\D$
and a certain quotient $\bV_\D$ of it. In a straightforward way one can compute
the elements of low degrees in the kernel of the projection
$$V_\D \rightarrow \bV_\D .$$
Applying the state-field correspondence map $Y$ to the elements of the kernel,
we obtain non-trivial relations on the action of the formal series
$$d_0 (m,z) = \sum_{j\in\ZS} t_0^j t_1^m d_0 z^{-j-2}$$
and
$$d_1 (m,z) = \sum_{j\in\ZS} t_0^j t_1^m d_1 z^{-j-1}.$$
In this way we derive the following relations in $\bV_\D$:
$$d_0(m,z) = 
{1\over 4} \lcol d_1(1,z) d_1(m-1,z)\rcol
- {1\over 2} \lcol d_1(0,z) d_1(m,z)\rcol
+ {1\over 4} \lcol d_1(-1,z) d_1(m+1,z)\rcol, \eqno{(\odz)}$$
$$\displaylines{
 {\d \over \d z} d_1(m,z)  =  
 - {m+1\over 2} \lcol d_1(1,z) d_1(m-1,z)\rcol \hfill \cr
\hfill + m \lcol d_1(0,z) d_1(m,z)\rcol
- {m-1\over 2} \lcol d_1(-1,z) d_1(m+1,z)\rcol . {(\odo)} \cr} $$

 We prove that the same relations hold in the simple modules $L(\alpha, \beta, \gamma)$,
provided that $\beta = -{\alpha (\alpha+1) \over 2}$. The relation (\odz)
is a generalization of the Sugawara construction and implies
that the action of $d_0(m,z)$ may be expressed by means of the action of the subalgebra
of the horizontal vector fields on the torus,
$$\H = \C[t_0^{\pm 1}, t_1^{\pm 1}] d_1 .$$
Hence, the module $L(\alpha, \beta, \gamma)$ 
with $\beta = -{\alpha (\alpha+1) \over 2}$,
remains irreducible when restricted to the subalgebra $\H$. 

 Our second step is to show that we can restrict to the yet smaller loop subalgebra
$\L$ and still have an irreducible module. We achieve this with the help 
of the differential equation (\odo).

 The idea of our proof may be seen from a well-known formula for the solution
of an ordinary differential equation
$$ {\d \over \d z} x(z) = a(z) x(z), \eqno{(\ode)}$$
which expresses $x(z)$ as 
$$x(z) = \exp \left( \int\limits_0^z a(z) dz \right) x(0) . \eqno{(\sole)}$$
This formula shows that the solution $x(z)$ is ``built'' from $a(z)$ and
the initial value $x(0)$. 

 We transform the differential equation (\odo) into an infinite system of 
differential equations for 
$P(m,n,z) = d_1(m,z) v(n)$, $m,n\in\Z$ (see (\vdom)). These differential 
equations are considerably more complicated than (\ode): it is an infinite
system of equations; there is more than one term in the right hand side;
the components of the operator-valued series in the right hand side
do not commute; it involves normally ordered products; the series contains 
positive and negative powers of $z$.

These equations are highly non-commutative, even compared to the Fock space 
situation, where the commutator relations in the Heisenberg algebra are
much simpler, and most of its components actually commute.

Because of this we need to redefine the exponential function, and this is 
done using iterative integration. To illustrate our method, we point out
that the solution (\sole) may be alternatively written as a series
$$x(z) = \sum_{k=0}^\infty b_k(z),$$
where $b_0 (z) = x(0)$ and $b_{k+1} (z) = \int\limits_0^z a(z) b_k (z) dz$.

In order to view an infinite system of equations as a single differential 
equation, we treat $m$ and $n$ as formal variables and work over the polynomial
algebra $\C [m,n]$.

The equation on $P(m,n,z)$ that we get will be of the form (see (\szp)):
$$ {\d \over \d z} P(m,n,z) = 
\left( z^{-1} \Omega + A_+(z) \right) P(m,n,z), \eqno{(\oap)}$$
where the operator $\Omega$ acts in the following way:
$$\displaylines{
 \Omega P(m,n,z) = \cr
\hfill - {1\over 2} (m+1)(n-\alpha) P(m-1,n+1,z)
+   mn P(m,n,z) - {1\over 2} (m-1)(n+\alpha) P(m+1,n-1,z) ,  
\hfill {(\ovo)} \cr} $$
and $A_+(z)$ is a power series in non-negative powers of $z$ with operator 
coefficients that correspond to the action of the loop subalgebra $\L$.

Another difficulty that we have to overcome here is that the term 
$z^{-1} \Omega$ can't be integrated in Laurent series. To handle this 
problem, we introduce a special integration formalism.

 We find the solution of (\oap) as a series
$$P(m,n,z) = \sum_{k=0}^\infty P_k (m,n,z) , \eqno{(\opi)}$$
where the terms in the right hand side are given by the recurrence
relation
$$P_{k+1} (m,n,z) = z^{\Omega} \int z^{-\Omega} A_+(z) P_k (m,n,z) dz \eqno{(\opii)}$$
with the integral defined in the following way:
$$ z^{\Omega} \int z^{-\Omega} \left( \sum_{j} R_j z^j \right) dz
= \sum_{j} \left( (j+1)I - \Omega \right)^{-1} R_j z^{j+1} .$$
In order for this solution to be well-defined, we require the operators 
$\Omega - (j+1)I$ to be invertible whenever $R_j \neq 0$.

The formula (\ovo) does not manifest invertibility of these operators 
on the space of polynomials $\C [m,n]$. In order to find the answer 
to this question, we give a representation-theoretic interpretation 
for the operator $\Omega$. It turns out that $\Omega$ may be viewed
as a Casimir operator on a tensor product of two Verma modules for
$\sll_2(\C)$.

We fix a parameter $\alpha \in \C$ and consider the space of polynomials 
$\C[x]$ as a module over the Witt algebra $\Der \C [t, t^{-1}]$
with the action given by
$$ (t^m d) p(x) = (\alpha m - x) p(x+m) .$$
We denote this module by $\C_\alpha [x]$.
Restricting to the subalgebra $\sll_2 (\C) = 
\C t^{-1}d \oplus \C d \oplus \C td$, we get a representation of
$\sll_2 (\C)$ on $\C_\alpha [x]$. 
In this framework, the operator $\Omega$ is just the bilinear 
Casimir operator on the tensor product 
$\C_{\alpha_1} [m] \otimes \C_{\alpha_2} [n]$ with 
$\alpha_1 = 1$, $\alpha_2 = \alpha$. 
A related construction of the action
of $\gll_n (\C)$ on the space of meromorphic functions in 
${1\over 2} n (n-1)$ variables, has been given in [GKL]. 

In spite of the fact that the action of degree zero derivation $d$
is not diagonalizable, $d p(x) = - x p(x)$, the module $\C_\alpha [x]$
is nonetheless isomorphic to a Verma module for $\sll_2 (\C)$
under a different choice of the Cartan element. This identification 
allows us to easily calculate the spectrum of $\Omega$
on the space $\C[m,n]$ and determine that for $\alpha\not\in\Q$
the operators $\Omega - jI$ are invertible for all $j = 0,1,2, \ldots$

 As a result, we conclude that the formulas (\opi) and (\opii) are valid,
and the components of the series $P(m,n,z)$ belong to the space
$U(\L_-) T$. Using this fact, it is easy to establish the isomorphism
(\isot).

The structure of the modules $L(\alpha,\beta,\gamma)$ for the values
of the parameters not covered by Theorem 1.2, remains an open problem. 

\

{\bf Acknowledgments.} The first author thanks the University 
of Sydney for hospitality during his visit there. 
Financial support from the Natural Sciences and
Engineering Research Council of Canada and the Australian Research Council
is gratefully acknowledged.

\

\

{\bf 1. Representations of Lie algebra of vector fields on a torus.}

\

 The algebra of Fourier polynomial functions on a 2-dimensional torus
is isomorphic to the algebra of Laurent polynomials in two variables
$\R = \C[t_0^{\pm 1}, t_1^{\pm 1}]$. The Lie algebra of (polynomial) vector 
fields on a torus is the algebra of derivations of $\R$:
$$ \D = \Der (\R) = \R {\d \over \d t_0} \oplus \R {\d \over \d t_1} .$$
It is convenient to use slightly different generators of $\Der (\R)$ as a
free $\R$-module:
$$ d_0 = t_0  {\d \over \d t_0}, \quad d_1 = t_1  {\d \over \d t_1} .$$

 The Lie bracket is given by the formula:
$$ \left[ t_0^{m_0} t_1^{m_1} d_a, t_0^{r_0} t_1^{r_1} d_b \right]
= r_a t_0^{m_0+r_0} t_1^{m_1+r_1} d_b - m_b t_0^{m_0+r_0} t_1^{m_1+r_1} d_a . \eqno{(\Lieb)}$$ 

 The degree operators $d_0, d_1$ induce a $\Z^2$-grading on $\D$ by their
eigenvalues in the adjoint representation. We will also consider a 
$\Z$-grading on $\D$ induced just by $d_0$.

 We will need two subalgebras in $\D$ -- the Lie algebra of ``horizontal''
vector fields
$$\H = \R d_1 $$
and the loop algebra
$$\L = \C[t_0, t_0^{-1}] \left\{ \C t_1^{-1} d_1 \oplus \C d_1 
\oplus \C t_1 d_1 \right\} .$$
It is easy to see that the latter subalgebra is isomorphic to the loop
Lie algebra
$$ \L \cong \C[t_0, t_0^{-1}] \otimes \sll_2 (\C) .$$

 Next we will construct a class of bounded modules for these Lie algebras.

 We take a decomposition of $\D$ into three parts,
$$\D = \D_+ \oplus \D_0 \oplus \D_- ,$$
corresponding to positive, zero and negative eigenvalues relative to
$d_0$. In particular,
$$\D_0 = \C[t_1, t_1^{-1}] d_0 \oplus \C[t_1, t_1^{-1}] d_1 .$$
In the same way we consider the decompositions for $\H$ and $\L$, 
where $\H_0 = \H \cap \D_0$, $\L_0 = \L \cap \D_0$, etc.
Note that $\H_0$ is the algebra of vector fields on a circle
and $\L_0 \cong \sll_2(\C)$.

 Consider the family of modules $T = T(\alpha, \beta, \gamma)$,
$\alpha,\beta,\gamma \in \C$, 
for $\D_0$ with the basis $\{v(n) | n\in \gamma+\Z \}$ and the action
defined by the formulas:
$$\eqalign{
(t_1^m d_1) v(n) &= (n - \alpha m) v(n+m) , \cr
(t_1^m d_0) v(n) &= \beta v(n+m). \cr} $$
 It is well-known that $T(\alpha, \beta, \gamma)$ is irreducible 
as a $\D_0$-module unless
$\alpha \in \{ -1, 0 \}, \beta = 0$ and $\gamma \in \Z$ (see e.g.~[E1]).

 We let $\D_+$ act on $T$ trivially, and construct
the generalized Verma module
$$ M_\D (T) = M_\D (\alpha, \beta, \gamma) = \Ind_{\D_0 \oplus \D_+}^\D 
T(\alpha, \beta,\gamma) \cong U(\D_-) \otimes T(\alpha, \beta,\gamma). $$
We call $T$ the {\it top} of $M_\D (T)$.

 In the same way, viewing $T$ as a module for $\H_0$ and $\L_0$,
we construct the generalized Verma modules over $\H$ and $\L$:
$$ M_\H (T) = M_\H (\alpha, \gamma) = 
U(\H_-) \otimes T, \quad
M_\L (T) = M_\L(\alpha, \gamma) = U(\L_-) \otimes T. $$
 In the notations above we dropped $\beta$ since the actions of $\H_0$
and $\L_0$ on $T(\alpha,\beta,\gamma)$ are independent of $\beta$. 

 These modules have $\Z^2$-gradings compatible with the gradings on the Lie algebras. 
In these gradings, the modules $M_\D (T)$ and $M_\H (T)$ have
infinite-dimensional components below the top, whereas all components of $M_\L (T)$ are
finite-dimensional.  

 The module $M_\D (T)$ has a unique maximal submodule. Indeed, 
any submodule in it is homogeneous
with respect to $\Z^2$-grading. Since $T$ as a $\D_0$-module has a 
unique maximal submodule (which is trivial in most cases), 
the intersection of any proper $\D$-submodule with
$T$ should be in the maximal submodule of $T$. Hence the sum of 
proper submodules of $M_\D (T)$ is again a proper submodule, and the sum
of all proper submodules is the unique maximal submodule.

The same argument shows that $M_\H (T)$ has a unique maximal homogeneous
(with respect to $\Z^2$-grading) submodule.
%
%
We define $L_\D (T) = L_\D (\alpha, \beta, \gamma)$ (resp. $L_\H (T)$) as a quotient
of $M_\D (T)$ (resp. $M_\H (T)$) by the maximal (resp. maximal homogeneous) submodule.


 The following Theorem is a direct corollary of a general result of
Berman and Billig ([BB], Theorem 1.12):

{\bf Theorem 1.1.} 
{\it
The modules $L_\D (\alpha, \beta, \gamma)$ and
$L_\H (\alpha, \gamma)$ have finite-dimensional components in
$\Z^2$-grading. 
}

\

 A natural question stems from Theorem 1.1:

{\bf Question.} 
{\it
What are the characters of the modules
$L_\D (\alpha, \beta, \gamma)$ and $L_\H (\alpha, \gamma)$?
}

 This question is also motivated by a conjecture which is due
to Eswara Rao [E2], which in the case of a 2-torus states: 

{\bf Conjecture.} 
{\it
Every simple $\Z^2$-graded module for the Lie algebra of 
vector fields on a $2$-dimensional torus with finite-dimensional
components of the grading is either a sub-quotient of a tensor module
or is isomorphic to $L_\D (\alpha, \beta, \gamma)$, possibly twisted
with an automorphism of $\D$.
}

 This conjecture is an analogue of Kac's conjecture, solved by Mathieu [M], which describes the case of the vector fields on a circle.

 Whereas the tensor modules are completely understood, nothing was previously 
known about the modules $L_\D (\alpha, \beta, \gamma)$. 
The goal of the present paper is to give an explicit description for 
some of these modules. We will prove the following

{\bf Theorem 1.2.} 
{\it
Let $\alpha \in \C, \alpha \not\in \Q$, 
$\beta = -{1 \over 2} \alpha (\alpha+1)$, $\gamma\in \C$. 
Then 
$$L_\D (\alpha, \beta, \gamma) \cong 
L_\H (\alpha, \gamma) \cong M_\L (\alpha, \gamma) .$$ 
That is, the action of the loop algebra $\L$ on $M_\L (\alpha, \gamma)$
can be extended to an irreducible action of the whole algebra of vector fields
$\D$ turning $M_\L (\alpha, \gamma)$ into 
$L_\D (\alpha, \beta, \gamma)$.
}

\


 Let $V$ be a module over the Lie algebra $\D$, graded by eigenvalues of
$d_0$ and $d_1$:
$$V = \mathop\bigoplus_{(s_0, s_1) \in \CS \times \CS} 
V_{s_0,s_1},$$
where 
$$ V_{s_0,s_1} =
\left\{ v \in V \, | \, d_0 v = s_0 v, d_1 v = s_1 v \right\},$$
and define the character of $V$ to be
$${\rm char\,} V = \sum\limits_{(s_0, s_1) \in \CS \times \CS} 
\dim \left( V_{s_0,s_1} \right) q_0^{s_0} q_1^{s_1} .$$

{\bf Corollary 1.3.} 
{\it
Under the assumptions of Theorem 1.2 on $\alpha$,
$\beta$, $\gamma$, 
$$ {\rm char\,} L_\D (\alpha, \beta, \gamma) = q_0^\beta  
\left( \sum_{j \in \gamma + \ZS} q_1^j \right) \prod\limits_{n\geq 1} (1-q_0^{-n})^{-3} . \eqno{(\charf)}$$
}

\

\

{\bf 2. Vertex algebras and vertex Lie algebras.}

\

 In this section we are going to construct vertex algebras associated with 
the Lie algebras $\D$, $\H$ and $\L$. 
We briefly recall some properties of vertex algebras and
refer to [K2] and [L] for details.

 A vertex algebra $V$ is a vector space with additional structures $(Y, D, \o)$, where
the state-field correspondence $Y$ is a map
$$ Y: V \rightarrow \End(V) [[z,z^{-1}]], $$
the infinitesimal translation $D$ is an operator $D: V \rightarrow V$,
and $\o$ is a vacuum vector $\o \in V$.

 For $a \in V$ we write
$$ Y(a,z) = \sum\limits_{n\in\ZS} a_{(n)} z^{-n-1}, \quad
a_{(n)} \in \End (V) .$$
For all $a,b \in V$,
$$ a_{(n)} b = 0 \quad {\rm for \ } n \gg 0.$$
The infinitesimal translation map $D$ satisfies the axiom
$$Y(Da, z) = {\d \over \d z} Y(a,z) .$$

Finally, we will need the commutator formula:
$$ [Y(a, z_1), Y(b, z_2) ] = \sum\limits_{n \geq 0} 
{1 \over n!} Y(a_{(n)} b, z_2)
\left[ z_1^{-1} \left( {\d \over \d z_2} \right)^n \delta \left(
{z_2 \over z_1} \right) \right]. \eqno{(\comm)}$$

 The delta function that appears in (\comm) above is
$$ \delta(z) = \sum_{n\in\ZS} z^n .$$

There exists a uniform procedure of constructing
vertex algebras from vertex Lie algebras, and the Lie algebras under consideration
are all vertex Lie algebras. Let us recall a definition of a vertex Lie algebra [DLM].

Let $\S$ be a Lie algebra with the basis 
$\{ u(n), c(-1) \big| u\in\U, c\in\CC, n\in\Z \}$ ($\U$, $\CC$ are some index sets).
Define the corresponding fields in $\S [[z,z^{-1}]]$:
$$ u(z) = \sum_{n\in\ZS} u(n) z^{-n-1}, \quad c(z) = c(-1) z^0, \quad 
u\in\U, c\in\CC .$$
Let $\F$ be a subspace in $\S [[z,z^{-1}]]$ spanned by all the fields
$u(z), c(z)$ and their derivatives of all orders.

{\bf Definition.}
{\it
 A Lie algebra $\S$ with the basis as above is called a vertex Lie algebra
if the following two conditions hold:

(VL1) for all $u_1, u_2 \in \U$,
$$ [u_1(z_1), u_2(z_2) ] = \sum\limits_{j=0}^n f_j(z_2)
\left[ z_1^{-1} \left( {\d \over \d z_2} \right)^j \delta \left(
{z_2 \over z_1} \right) \right], \eqno{(\vla)}$$
where $f_j(z) \in\F, n \geq 0$ and both depend on $u_1, u_2$,

(VL2) for all $c\in\CC$, the elements $c(-1)$ are central in $\S$.
}

\

Let $\S^{(+)}$ be a subspace in $\S$ with the basis $\{ u(n) \big| u\in\U, n\geq 0 \}$
and let  $\S^{(-)}$ be a subspace with the basis 
$\{ u(n), c(-1) \big| u\in\U, c\in\CC, n<0 \}$. Then $\S = \S^{(+)} \oplus \S^{(-)}$ and
$\S^{(+)}, \S^{(-)}$ are in fact subalgebras in $\S$.

The universal enveloping vertex algebra $V_{\S}$ of a vertex Lie algebra $\S$ 
is defined as an induced module
$$V_{\S} = \Ind_{\S^{(+)}}^\S (\C \o) = U(\S^{(-)}) \ot \o,$$
where $\C \o$ is a trivial 1-dimensional $\S^{(+)}$ module.

{\bf Theorem \voa. ([DLM], Theorem 4.8)} 
{\it
Let $\S$ be a vertex Lie algebra. Then

(a) $V_{\S}$ has a structure of a vertex algebra with the vacuum vector $\o$,
infinitesimal translation $D$ being a natural extension of 
the derivation of $\S$ given by 
$$D(u(n)) = -n u(n-1),  \quad D(c(-1)) = 0, \quad u\in\U, c\in\CC,$$
and the state-field correspondence map $Y$ defined by the formula:
$$Y \big( a_1(-1-n_1) \ldots a_{k-1}(-1-n_{k-1}) a_k(-1-n_k) \o, z
\big) $$
$$ = 
\lcol \left( {1\over n_1 !} \left( {\d \over \d z} \right)^{n_1} a_1 (z)
\right) 
\ldots
\lcol \left( {1\over n_{k-1} !} \left( {\d \over \d z} \right)^{n_{k-1}} 
a_{k-1} (z) \right) 
\left( {1\over n_{k} !} \left( {\d \over \d z} \right)^{n_k} 
a_k (z) \right)\rcol \ldots \rcol \quad 
, \eqno{(\Y)}$$
where $a_j \in \U, n_j \geq 0$ or $a_j \in\CC, n_j =0$.

(b) Any bounded $\S$-module is a vertex algebra module for $V_{\S}$.
}

\

In the formula (\Y) above,  
the normal
ordering of two fields $\lcol a(z) b(z) \rcol$ is defined as 
$$ \lcol a(z) b(z) \rcol = a_+(z) b(z) + b(z) a_-(z) ,$$
where $a_+(z)$ (resp. $a_-(z)$) is the part of the power series $a(z)$
involving non-negative (resp. negative) powers of $z$,
$a(z) = a_+(z) + a_-(z)$.

\

 There is a natural relation between $\S$-submodules of $V_\S$ and vertex algebra
ideals of $V_\S$.

{\bf Theorem \inv. (see e.g.~[B])} 
{\it Any $D$-invariant $\S$-submodule in $V_{\S}$ is a vertex algebra
ideal in $V_{\S}$. Conversely, every vertex algebra ideal in $V_{\S}$ is a 
$D$-invariant $\S$-submodule.

}
 
\

{\bf Theorem \vl.} 
{\it 
The Lie algebras $\D$, $\H$ and $\L$ are vertex Lie algebras.
The generating fields in $\D$ are
$$ d_1 (m,z) = \sum_{j\in \ZS} t_0^j t_1^m d_1 z^{-j-1}, \eqno{(\done)} $$
$$ d_0 (m,z) = \sum_{j\in \ZS} t_0^j t_1^m d_0 z^{-j-2}, \eqno{(\dzer)} $$
where $m\in \Z$.
The generating fields of $\H$ are (\done) with $m\in \Z$, and the generating fields
of $\L$ are (\done) with $m = -1, 0, 1$. 
}

Note that the centers of these Lie algebras are trivial, and the set $\CC$ is empty
in all of these cases. 

{\bf Proof.}
We need to show that the commutators between the fields (\done) and (\dzer) conform
with (\vla).

$$ \left[ d_1 (m,z_1), d_1(r, z_2) \right] = 
\sum_{ij} \left[ t_0^i t_1^m d_1, t_0^j t_1^r d_1 \right] z_1^{-i-1} z_2^{-j-1} $$
$$ = (r-m) \sum_{ij} \left( t_0^{i+j} t_1^{m+r} d_1 z_2^{-i-j-1} \right) 
\left( z_1^{-i-1} z_2^i \right) = 
(r-m) d_1(m+r,z_2) \z . \eqno{(\doomr)}$$

$$ \left[ d_1 (m,z_1), d_0(r, z_2) \right] = 
\sum_{ij} \left[ t_0^i t_1^m d_1, t_0^j t_1^r d_0 \right] z_1^{-i-1} z_2^{-j-2} $$
$$ = r \sum_{ij} \left( t_0^{i+j} t_1^{m+r} d_0 z_2^{-i-j-2} \right) 
\left( z_1^{-i-1} z_2^i \right)
 - \sum_{ij} \left( t_0^{i+j} t_1^{m+r} d_1 z_2^{-i-j-1} \right) 
\left( i z_1^{-i-1} z_2^{i-1} \right)$$
$$ = r d_0 (m+r, z_2) \z - d_1 (m+r, z_2) \zd .  \eqno{(\dozmr)}$$

$$ \left[ d_0 (m,z_1), d_0(r, z_2) \right] = 
\sum_{ij} \left[ t_0^i t_1^m d_0, t_0^j t_1^r d_0 \right] z_1^{-i-2} z_2^{-j-2} $$
$$ =  \sum_{ij} (j-i) t_0^{i+j} t_1^{m+r} d_0 z_1^{-i-2} z_2^{-j-2} $$
$$ = \sum_{ij} \left( (i+j+2) t_0^{i+j} t_1^{m+r} d_0 z_2^{-i-j-3} \right) 
\left( z_1^{-i-2} z_2^{i+1} \right)$$
$$ - 2 \sum_{ij} \left( t_0^{i+j} t_1^{m+r} d_0 z_2^{-i-j-2} \right) 
\left( (i+1) z_1^{-i-2} z_2^i \right)$$
$$ = - \left( {\d \over \d z_2} d_0 (m+r, z_2) \right) \z - 2 d_0 (m+r, z_2) \zd .  \eqno{(\dzzmr)}$$

This proves the claim of the theorem.
  
 \

It follows from (\done), (\dzer) that
$$ \D^{(+)} = \Span \left\{ t_0^i t_1^m d_1, t_0^j t_1^m d_0 \; | \;
i \geq 0, j \geq -1, m\in \Z \right\},$$
$$ \D^{(-)} = \Span \left\{ t_0^i t_1^m d_1, t_0^j t_1^m d_0 \; | \;
i \leq -1, j \leq -2, m\in \Z \right\},$$
and $\H^{(\pm)} = \H \cap \D^{(\pm)}$, $\L^{(\pm)} = \L \cap \D^{(\pm)}$.

 By applying Theorem \voa, we obtain the vertex algebras 
$V_\D = U(\D^{(-)}) \otimes \o$,
$V_\H = U(\H^{(-)}) \otimes \o$ and $V_\L = U(\L^{(-)}) \otimes \o$. 
We have natural embeddings of these vertex algebras:
$V_\L \subset V_\H \subset V_\D$.

\

\

{\bf 3. Relations in vertex algebras and modules.}

\

 The vertex algebras $V_\D$, $V_\H$ and $V_\L$ are graded by $\Z^2$ with $\o$ having 
degree $(0,0)$. In the $\Z$-grading by the degree in $t_0$, these algebras are trivial in 
positive degree. At the top, in degree zero relative to $t_0$, we have
the only non-trivial component $\C\o$. When we go to the next layer
in degree $-1$, then in $V_\D$ and $V_\H$ we get an infinite direct sum
of 1-dimensional spaces spanned by $(t_0^{-1} t_1^m d_1) \o$ in each degree $(-1,m)$. 
Note the following relation in $V_\D$:
$$ (t_0^{-1} t_1^m d_0) \o = 0 .\eqno{(\doo)} $$

 The vertex algebras $V_\D$, $V_\H$ and $V_\L$ are not simple. In fact, their maximal 
ideals correspond to the augmentation ideals in
$U(\D^{(-)})$, etc., and thus the quotients by the maximal ideals yield
just a 1-dimensional trivial vertex algebra $\C \o$. Nonetheless, 
it is possible to extract useful 
information for the representation theory of the
corresponding Lie algebras by studying certain quotients of
these vertex algebras.

 Let $V$ be one of the vertex algebras, $V_\D$, $V_\H$ or $V_\L$.
As an alternative to the maximal ideal, we will define the 
{\it subradical} of $V$ as an ideal $V^{sr}$ in $V$, maximal among the
homogeneous (in $\Z^2$-grading)
ideals trivially intersecting with the space $V_{-1}$ of elements of degree $-1$ with respect 
to $t_0$. 
 
 Let $\bV$ be the quotient $ \bV = V / V^{sr} .$

  {\bf Theorem 3.1.} 
{\it
The vertex algebras $\bV_\D$ and $\bV_\H$ are isomorphic. More precisely, 
there exists an isomorphism
$\bV_\H \rightarrow \bV_\D$, making the following diagram commutative:
}


$$\matrix{ V_\H & \longrightarrow & V_\D \cr
\Big\downarrow && \Big\downarrow \cr
\bV_\H & \longrightarrow & \bV_\D \cr }$$

 This theorem implies that the $\H$-module $\bV_\H$ admits a natural 
structure of a module over the larger Lie algebra $\D$.

 In order to prove this theorem, we need to study the subradical
$V_\D^{sr}$ at degree $-2$. The elements of the subradical yield
relations in $\bV_\D$ which will be crucial for our work.

 {\bf Proposition 3.2.} 
{\it
(a) The component of $\bV_\D$ of degree $(-2,m)$ has dimension 3, and is spanned 
by the elements $(t_0^{-1} t_1^k d_1) (t_0^{-1} t_1^{m-k} d_1) \o$ with $k = -1, 0 ,1$.

(b) The following relations in $\bV_\D$ describe the projection
$V_\D \rightarrow \bV_\D$ in degree $(-2,m)$:
}
$$\displaylines{
 (t_0^{-2} t_1^m d_0) \o = \cr
\hfill
{1\over 4} (t_0^{-1} t_1 d_1) (t_0^{-1} t_1^{m-1} d_1) \o 
 - {1\over 2} (t_0^{-1} d_1) (t_0^{-1} t_1^{m} d_1) \o 
+ {1\over 4} (t_0^{-1} t_1^{-1} d_1) (t_0^{-1} t_1^{m+1} d_1) \o,
\hfill \quad{(\dzm)} \cr} $$
$$ \displaylines{
(t_0^{-2} t_1^m d_1) \o =  
- {m+1\over 2} (t_0^{-1} t_1 d_1) (t_0^{-1} t_1^{m-1} d_1) \o
\hfill \cr
\hfill + m (t_0^{-1} d_1) (t_0^{-1} t_1^{m} d_1) \o 
- {m-1\over 2} (t_0^{-1} t_1^{-1} d_1) (t_0^{-1} t_1^{m+1} d_1) \o,
\quad{(\dom)} \cr} $$
$$\displaylines{
 (t_0^{-1} t_1^r d_1) (t_0^{-1} t_1^{m-r} d_1) \o = 
 {r(r+1)\over 2} (t_0^{-1} t_1 d_1) (t_0^{-1} t_1^{m-1} d_1) \o 
\hfill\cr
\hfill - (r-1)(r+1) (t_0^{-1} d_1) (t_0^{-1} t_1^{m} d_1) \o 
+{r(r-1) \over 2} (t_0^{-1} t_1^{-1} d_1) (t_0^{-1} t_1^{m+1} d_1) \o.
\quad{(\drm)} \cr} $$  

{\bf Proof.}
In order to prove this proposition we need to study the action of the raising operators 
$t_0 t_1^s d_0$ and $t_0 t_1^s d_1$ on the component 
of degree $(-2,m)$ of $V_\D$. The corresponding component of the
subradical $V_\D^{sr}$ is the joint kernel of these raising operators.

 The following calculations are straightforward and use (\doo):
$$ \eqalign{
(t_0 t_1^s d_0) (t_0^{-2} t_1^m d_0) \o &= 0, \cr
(t_0 t_1^s d_1) (t_0^{-2} t_1^m d_0) \o &= - (t_0^{-1} t_1^{m+s} d_1) \o,} $$
$$ \eqalign{
(t_0 t_1^s d_0) (t_0^{-2} t_1^m d_1) \o &= - 2(t_0^{-1} t_1^{m+s} d_1) \o, \cr
(t_0 t_1^s d_1) (t_0^{-2} t_1^m d_1) \o &= (m-s) (t_0^{-1} t_1^{m+s} d_1) \o,} $$
$$\displaylines{
\hfill
(t_0 t_1^s d_0) (t_0^{-1} t_1^r d_1) (t_0^{-1} t_1^{m-r} d_1) \o = 
(2s+2r-m) (t_0^{-1} t_1^{m+s} d_1) \o, 
\hfill
\hbox{\hskip 1.5cm}
\cr
\hfill
(t_0 t_1^s d_1) (t_0^{-1} t_1^r d_1) (t_0^{-1} t_1^{m-r} d_1) \o = 
(r-s)(m-2r-s) (t_0^{-1} t_1^{m+s} d_1) \o. 
 \hfill (\rais) \cr}  $$

 Since $s$ is an arbitrary integer, we will treat the coefficients in the
right hand sides above as polynomials in $s$. 
We consider a map
$\psi$ from the $(-2, m)$-component of $V_\D$ to $\C[s] \oplus \C[s]$
where the two components of $\psi (x)$ are the coefficients at 
$(t_0^{-1} t_1^{m+s} d_1) \o$ in $(t_0 t_1^s d_0) x$ and
$(t_0 t_1^s d_1) x$. Then the kernel of $\psi$ consists of the elements in 
the subradical, and the dimension of the quotient is equal to the
dimension of the image of $\psi$. 

 By inspection, we can see that the image of $\psi$ is spanned by
$\{(0,1), (2,s), (2s+m, s^2) \}$, thus the quotient by the subradical
in degree $(-2,m)$ has dimension 3. It is easy to verify that
$\psi((t_0^{-1} t_1^r d_1) (t_0^{-1} t_1^{m-r} d_1) \o)$ with
$r=-1,0,1$ span the image of $\psi$. Using the formulas above it is easy 
to check that the differences between left hand sides and right hand sides 
in (\dzm)-(\drm) are in the subradical.
 
\

 {\bf Remark.} In fact (\dom) can be derived from (\drm) by using the 
relation
$t_0^{-2} t_1^m d_1 = (m-2r)^{-1} 
[t_0^{-1} t_1^r d_1, t_0^{-1} t_1^{m-r} d_1]$
with $m-2r \neq 0$.

 \

{\bf Corollary 3.3.}
{\it
The following relations hold in the vertex algebra $\bV_\D$:}
$$d_0(m,z) = 
{1\over 4} \lcol d_1(1,z) d_1(m-1,z)\rcol
- {1\over 2} \lcol d_1(0,z) d_1(m,z)\rcol 
+ {1\over 4} \lcol d_1(-1,z) d_1(m+1,z)\rcol,
\eqno{(\Ydzm)} $$
$$\displaylines{
 {\d \over \d z} d_1(m,z)  =  
 - {m+1\over 2} \lcol d_1(1,z) d_1(m-1,z)\rcol \hfill \cr
\hfill + m \lcol d_1(0,z) d_1(m,z)\rcol 
- {m-1\over 2} \lcol d_1(-1,z) d_1(m+1,z)\rcol , 
\quad{(\Ydom)} \cr} $$
$$\displaylines{
 \lcol d_1(r,z) d_1(m-r,z)\rcol =  {r(r+1)\over 2} \lcol d_1(1,z) d_1(m-1,z)\rcol \hfill \cr
\hfill - (r-1)(r+1) \lcol d_1(0,z) d_1(m,z)\rcol 
+{r(r-1) \over 2} \lcol d_1(-1,z) d_1(m+1,z)\rcol .
\quad{(\Ydrm)} \cr} $$  

{\bf Proof.} The relations (\Ydzm)-(\Ydrm) are derived from (\dzm)-(\drm)
by applying the state-field correspondence map $Y$, and noting that
from the construction of the universal enveloping vertex algebra we have
$Y((t_0^{-2} t_1^m d_0)\o, z) = d_0(m,z)$,
$Y((t_0^{-1} t_1^m d_1)\o, z) = d_1(m,z)$ and
$Y((t_0^{-2} t_1^m d_1)\o, z) = {\d \over \d z} d_1(m,z)$
both in $V_\D$ and $\bV_\D$.

\

Now we are ready to give a proof of Theorem 3.1. It follows from (\Ydzm)
that for any $v \in \bV_\D$ we have $(t_0^j t_1^m d_0) v \in
U(\H) v$ for all $j,m \in \Z$. This implies that the map
$$V_\H \rightarrow \bV_\D ,$$
obtained as a composition of embedding $V_\H \rightarrow V_\D$
and projection $V_\D \rightarrow \bV_\D$, is surjective. It remains 
to show that the kernel of this map is the subradical of $V_\H$,
which follows from 

{\bf Lemma 3.4.} $V_\D^{sr} \cap V_\H = V_\H^{sr}$.

{\bf Proof.} Indeed, let $v$ be a homogeneous element in 
$V_\H \subset V_\D$. This element belongs to
$V_\D^{sr}$ if and only if $U(\D) v$ trivially intersects with the
degree $-1$ layer of $V_\D$. Note that the degree $-1$ layers in 
$V_\D$ and $V_\H$ coincide.
Due to (\Ydzm), $U(\D) v = U(\H) v \; 
\mod V_\D^{sr}$. Since $V_\D^{sr}$ trivially intersects with the
degree $-1$ layer of $V_\D$, we get that $v \in V_\D^{sr}$ if and 
only if $v \in V_\H^{sr}$. This completes the proof of the Lemma.

\

Thus the map $V_\H \rightarrow \bV_\D$
factors through $\bV_\H$ and the resulting map 
$$\bV_\H \rightarrow \bV_\D$$
is an isomorphism of $\H$-modules. Theorem 3.1 is now proved.

\

We will prove below an analogue of Theorem 3.1 for the modules over
$\D$ and $\H$. We will show that irreducible $\H$-modules $L_\H(\alpha,\gamma)$
admit the action of the bigger algebra $\D$, with the action extended 
by (\Ydzm). In order to prove this result, we will need a description
of the layer of degree $-3$ in $\bV_\H$. Since the subalgebra $\H^{(-)}$ is generated
by its elements of degree $-1$, the component of degree $(-3,m)$ in $V_\H$ is
spanned by the elements $(t_0^{-1} t_1^r d_1) (t_0^{-1} t_1^k d_1) (t_0^{-1} t_1^{m-r-k} d_1) \o$,
and  hence we will need to describe how such elements reduce in $\bV_\H$.

\

 {\bf Proposition 3.5.} 
{\it
(a) The component of $\bV_\H$ of degree $(-3,m)$ has dimension 8, and is spanned 
by the elements $(t_0^{-1} t_1^i d_1) (t_0^{-1} t_1^j d_1) (t_0^{-1} t_1^{m-i-j} d_1) \o$ 
with $i,j = -1, 0 ,1$, which satisfy the following linear relation in $\bV_\H$:
}
$$ \sum_{i,j=-1,0,1} A_{i+2, j+2} (t_0^{-1} t_1^i d_1) (t_0^{-1} t_1^j d_1) (t_0^{-1} t_1^{m-i-j} d_1) \o = 0, \eqno{(\relm)}$$
{\it
with
}
$$A = \pmatrix{
m(m-1)(m-2) & -2m(m-1)(m-2) & (m-3)(m-2)(m+2) \cr
-2(m-1)(m-2)(m+3) & 4m (m-2)(m+2) & -2(m-3)(m+1)(m+2) \cr
(m-2)(m+2)(m+3) & -2m(m+1)(m+2) & m(m+1)(m+2) \cr
} . $$ 


{\it
(b) The following relation holds in $\bV_\D$ and describes the projection
$V_\H \rightarrow \bV_\H$ in degree $(-3,m)$:
}
$$ (t_0^{-1} t_1^r d_1) (t_0^{-1} t_1^k d_1) (t_0^{-1} t_1^{m-r-k} d_1) \o = 
\sum_{i,j=-1,0,1} B_{i+2, j+2} (t_0^{-1} t_1^i d_1) (t_0^{-1} t_1^j d_1) (t_0^{-1} t_1^{m-i-j} d_1) \o \eqno{(\drkm)}$$
{\it
with
}
$$\displaylines{
B = 
{1\over 12} \pmatrix{
r(r-1) & 0 & 0 \cr
0 & (r-1)(r+1) & 0 \cr
0 & 0 & r(r+1) \cr} \times \hfill \cr
\hfill \times \pmatrix{
c+r-k & -2(c+4r+2k) & c+7r+5k+6 \cr
-2(c-3r-3k) & 4(c-3) & -2(c+3r+3k) \cr
c-7r-5k+6 & -2(c-4r-2k) & c-r+k \cr
}, \cr}$$    
%
where $c = r^2 + 2rk + 3k^2$.

\

 This proposition is proved in the same way as Proposition 3.2, and we will omit the details of this
straightforward calculation.

 Let $R$ be the vertex algebra ideal in $V_\H$ generated by the elements of $V_\H^{sr}$ corresponding to the relations
(\drm), (\relm) and (\drkm). Define $\tV_\H$ to be the quotient $V_\H /R$. The reason for considering such a quotient
is that we need relations corresponding to the elements of degrees $-2$ and $-3$ in $V^{sr}_\H$, and yet 
when we study modules for $V_\H$, it is easier to check that the relations corresponding to (\drm), (\relm) and (\drkm)
hold in the given module, rather than checking that all relations corresponding to the subradical $V^{sr}_\H$ hold.

\

{\bf Theorem 3.6.} 
{\it
There is an epimorphism of vertex algebras
$$ V_\D \rightarrow \tV_\H$$
which is defined by (\dzm) on the generators $(t_0^{-2} t_1^m d_0) \o$ 
and by the identity map on 
\break
$(t_0^{-1} t_1^m d_1) \o$.
}

\

{\bf Proof.} Denote the right hand side of (\dzm) by $\td_0(m)$. We need to show that the fields
$Y(\td_0 (m), z)$ and $Y( (t_0^{-1} t_1^m d_1) \o, z)$ in $\tV_\H$ yield a representation of the Lie algebra $\D$,
i.e., they satisfy relations (\doomr)-(\dzzmr). 
We can use the commutator formula (\comm) to express these relations
via $n$-th products:
$$ (t_0^{-1} t_1^m d_1 \o)_{(0)} (t_0^{-1} t_1^r d_1 \o) = (r-m) (t_0^{-1} t_1^{r+m} d_1 \o) , \eqno{(\dooz)}$$
$$ (t_0^{-1} t_1^m d_1 \o)_{(n)} (t_0^{-1} t_1^r d_1 \o) = 0 \quad \hbox{\rm \ for \ } n\geq 1 , \eqno{(\dooo)}$$
$$ (t_0^{-1} t_1^m d_1 \o)_{(0)} \td_0(r) = r \td_0 (r+m) , \eqno{(\dozz)}$$
$$ (t_0^{-1} t_1^m d_1 \o)_{(1)} \td_0(r) = - (t_0^{-1} t_1^{r+m} d_1 \o) , \eqno{(\dozo)}$$
$$ (t_0^{-1} t_1^m d_1 \o)_{(n)} \td_0(r) = 0 \quad \hbox{\rm \ for \ } n\geq 2 , \eqno{(\dozt)}$$
$$ \td_0(m)_{(0)} \td_0(r) = -D \td_0 (r+m) , \eqno{(\dzzz)}$$
$$ \td_0(m)_{(1)} \td_0(r) = -2 \td_0 (r+m) , \eqno{(\dzzo)}$$
$$ \td_0(m)_{(n)} \td_0(r) = 0 \quad \hbox{\rm \ for \ } n\geq 2 . \eqno{(\dzzt)}$$

Relations (\dooz) and (\dooo) obviously hold. To verify (\dzzz), we will work in the vertex algebra
$V_\D$. We have 
$$ \td_0(m)_{(0)} \td_0(r) =  (t_0^{-2} t_1^m d_0 \o)_{(0)} (t_0^{-2} t_1^r d_0 \o) \quad \mod V^{sr}_\D$$
$$ = -D (t_0^{-2} t_1^{r+m} d_0 \o) = -D \td_0(r+m) \quad \mod V^{sr}_\D.$$
Thus the difference between the right hand side and the left hand side in (\dzzz) belongs to $V^{sr}_\D$.
By Lemma 3.4 it is actually in $V^{sr}_\H$. Note also that these elements are of degree $-3$, and the ideal
$R$ contains all elements of $V^{sr}_\H$ of degree $-3$ by its construction. Hence the two sides of
(\dzzz) are equal in $\tV_\H$. The verification of other relations is completely analogous.

\

{\bf Corollary 3.7.} 
{\it
Every module for the vertex algebra $\tV_\H$ admits the action of the Lie algebra $\D$ defined by (\Ydzm).
}

\

{\bf Theorem 3.8.} 
{\it
Let $\beta = -{\alpha (\alpha + 1) \over 2}$. Then the modules $L_\H (\alpha, \gamma)$ and
$L_\D (\alpha, \beta, \gamma)$ are isomorphic. The action of $\H$ on $L_\H (\alpha, \gamma)$ extends to
the action of $\D$ by (\Ydzm).
}

 The idea of the proof is to show that $L_\H (\alpha, \gamma)$ is a module for the vertex algebra
$\tV_\H$ and then apply Corollary 3.7. To show that $L_\H (\alpha, \gamma)$ admits the action of
$\tV_\H$, we need to prove that for every generator $u$ of the ideal $R$ in $V_\H$, we have 
$$Y(u,z) L_\H (\alpha, \gamma) = 0.$$
Expand
$$ Y(u,z) = \sum\limits_{j\in \ZS} u_{(j)} z^{-j-1}, \quad u_{(j)} \in \End( L_\H (\alpha, \gamma) ) .$$
Consider the subspace $S$ in $\End( L_\H (\alpha, \gamma) )$ spanned by $u_{(j)}$, $j\in\Z$, where $u$
runs over the elements of $V^{sr}_\H$ of degrees $-2$ and $-3$.

\

{\bf Lemma 3.9.} 
{\it
The space $S$ is invariant with respect to the adjoint action of $\H$.
}

{\bf Proof.} The adjoint action of $\H$ may be expressed using $n$-th products:
$$ [t_0^i t_1^m d_1, u_{(j)} ] = [(t_0^{-1} t_1^m d_1 \o)_{(i)}, u_{(j)} ] .$$
To deal with the last expression we recall the Borcherds commutator formula:
$$[a_{(i)}, b_{(j)}] = \sum\limits_{k\geq 0} 
\pmatrix{ i \cr k \cr} (a_{(k)} b)_{(i+j-k)} .$$
Note that the operators $(t_0^{-1} t_1^m d_1 \o)_{(k)}$ with $k\geq 0$ increase the
degree by $k$. Since the subradical $V^{sr}_\H$ is a vertex algebra ideal, which is trivial 
in degree greater or equal to $-1$, we get that the span of elements of $V^{sr}_\H$
of degrees $-2$ and $-3$ is stable under the action of $(t_0^{-1} t_1^m d_1 \o)_{(k)}$ with $k\geq 0$.
Thus for $k \geq 0$,  $\left( (t_0^{-1} t_1^m d_1 \o)_{(k)} u \right)_{(i+j-k)}$ belongs to $S$
and the Lemma is proved.

\

{\bf Lemma 3.10.} 
{\it
The operators in $S$ of degree $0$ act trivially on the top $T$ of the module
$L_\H (\alpha, \gamma)$.  
}

{\bf Proof.} The space $S$ is spanned by the moments of the following generating series:
$$ \lcol d_1(r,z) d_1(m-r)\rcol 
- {r(r+1)\over 2} \lcol d_1(1,z) d_1(m-1,z)\rcol +$$
$$+ (r-1)(r+1) \lcol d_1(0,z) d_1(m,z)\rcol 
- {r(r-1) \over 2} \lcol d_1(-1,z) d_1(m+1,z)\rcol ,
\eqno{(\fdrm)} $$  
$$ \lcol d_1(r,z) \lcol d_1(k,z) d_1(m-r-k,z)\rcol\rcol 
- \sum_{i,j=-1,0,1} B_{i+2,j+2} \lcol d_1(i,z) \lcol d_1(j,z) d_1(m-i-j,z)\rcol\rcol , 
\eqno{(\frkm)} $$
%
and
$$\sum_{i,j=-1,0,1} A_{i+2,j+2} 
\lcol d_1(i,z) \lcol d_1(j,z) d_1(m-i-j,z)\rcol \rcol . \eqno{(\fmm)} $$
%


The degree zero component of $\lcol d_1 (p,z) d_1(s,z)\rcol$ is the $z^{-2}$ moment of this
formal series, and it acts on $v(n)$ as
$$ (t_1^s d_1) (t_1^p d_1) v(n) = (n-\alpha p)(n+p-\alpha s) v(n+p+s).$$
Verifying that the zero moments of (\fdrm), (\frkm)
and (\fmm) vanish on $T$, amounts to checking that the corresponding sums of polynomials are zero.
This can be easily done with {\it Maple} or any other software for symbolic computation,
or directly by hand if the reader favours the traditional method.
This completes the proof of the Lemma.

\ 

{\bf Proof of Theorem 3.8.} Let $u$ be a homogeneous element of the subradical 
$V^{sr}_\H$ of degree $-2$ or $-3$. We need to show that the operators $u_{(n)}$ vanish in 
$L_\H(\alpha, \gamma)$.
Let $\sum_k w_k v(k)$ be a homogeneous element in $L_\H (\alpha,\gamma)$, $w_k \in U(\H_-)$.
To prove that $u_{(n)} \sum_k w_k v(k) = 0$, we need to show that for every element $w^\prime \in U(\H_+)$
we get 
\break
$w^\prime u_{(n)} \sum_k w_k v(k) = 0$ whenever 
the left hand side belongs to the top $T$.

Using the Poincar\'e-Birkhoff-Witt argument, we can move all raising operators in the product
$w^\prime u_{(n)} w_k$ to the right, and all lowering operators to the left, keeping degree zero
operators in the middle. By Lemma 3.9, each term will contain a factor from $S$. 
 Applying such an expression to $v(k)$ we see that all the summands where the
raising operators are present, vanish. The terms that remain will be products of degree zero operators.
But Lemma 3.10 implies that every degree zero operator from $S$ vanishes on $T$. Thus
$w^\prime u_{(n)} \sum_k w_k v(k) = 0$ in $L_\H(\alpha, \gamma)$, and so $u_{(n)}$ vanishes in this module.
This proves that $L_\H(\alpha, \gamma)$ is a module for $\tV_\H$, and hence by Corollary 3.7,
the action of $\H$ on $L_\H(\alpha, \gamma)$ extends to the action of $\D$ on the same space using (\Ydzm).

In this way we get an irreducible module over $\D$. In order to identify it with one of the modules
$L_\D(\alpha,\beta,\gamma)$, we need to compute the action of $t_1^m d_0$ on $T$:
$$(t_1^m d_0) v(n) =  {1\over 4} (t_1^{m-1} d_1) (t_1 d_1) v(n) - {1\over 2}  (t_1^{m} d_1) d_1 v(n)
+ {1\over 4}  (t_1^{m+1} d_1) (t_1^{-1} d_1) v(n) $$
$$ = \left( {1\over 4}  (n-\alpha) (n+1 - \alpha(m-1)) -{1\over 2}  n (n - \alpha m) + {1\over 4} (n+\alpha) (n-1 - \alpha(m+1))
\right) v(n+m) $$
$$= - {\alpha (\alpha +1) \over 2} v(n+m) .$$
Thus $L_\H(\alpha, \gamma)$  is isomorphic to $L_\D(\alpha, \beta, \gamma)$ with $\beta = - {\alpha (\alpha +1) \over 2}$.
This completes the proof of the theorem. 

\

{\bf Corollary 3.11.} 
{\it
The relations (\Ydom) hold in $L_\H (\alpha, \gamma)$.
}

{\bf Proof.}
 In the proof of Theorem 3.8, we actually showed that all the relations corresponding to the elements of the subradical
$V_\H^{sr}$ of degree $-2$ and $-3$, hold in the modules  $L_\H (\alpha, \gamma)$. 

\

We will view (\Ydom) as a system
of differential equations. Solving these differential equations will allow us to prove that for $\alpha \not\in\Q$ the module $L_\H (\alpha, \gamma)$
is generated by $T$ as a module over the loop subalgebra $\L$.

\

\

{\bf 4. Representations of $\sll_2 (\C)$}.

\

In order to solve the differential equations (\Ydom), we need to study a family of representations for the Lie
algebra $\sll_2(\C)$. Fix $\alpha\in\C$, and consider the following action of the Lie algebra of vector fields
on a circle $\Der \C[t,t^{-1}]$ on the space of polynomials $\C[x]$:
$$ (t^m d) p(x) = (\alpha m - x) p(x+m) . \eqno{(\Wact)}$$
Let us verify that this is indeed a representation of the Lie algebra:
$$ (t^m d) (t^r d) p(x) - (t^r d) (t^m d) p(x) $$
$$= (\alpha m-x) (\alpha r -x -m) p(x+m+r) - (\alpha r-x) (\alpha m-x-r) p(x+m+r)$$
$$= (r-m) (\alpha(m+r)-x) p(x+m+r) = [t^m d, t^r d] p(x) .$$
 We will denote this module by $\C_\alpha [x]$.

 Let us restrict the module $\C_\alpha [x]$ to the subalgebra $\C t^{-1} d  \oplus \C d \oplus \C t d$, which is isomorphic to
$\sll_2 (\C)$. Note that the Cartan subalgebra element $d$ is not diagonalizable on $\C_\alpha[x]$, as it acts by 
$d p(x) = - x p(x)$. In particular, this shows that $\C_\alpha [x]$ is a cyclic module, generated by vector $1$ under the action
of $d$. Even though $d$ is not diagonalizable on $\C_\alpha [x]$, still this module is isomorphic to a Verma module over $\sll_2 (\C)$,
but with respect to a non-standard Cartan subalgebra.

 {\bf Theorem 4.1.} 
{\it
The module $\C_\alpha[x]$ is isomorphic to the Verma module for $\sll_2 (\C)$ with the highest weight $2\alpha$
and the highest weight vector $1$
with respect to the following Cartan decomposition of $\sll_2(\C)$: 
$$ h = td - t^{-1}d, \quad e = {1\over 2} \left( td - 2d + t^{-1} d \right), \quad
f = - {1\over 2} \left( td + 2d + t^{-1} d \right) .$$
}

{\bf Proof.} Verification of the usual relations $[h,e] = 2e$, $[h,f] = -2f$, $[e,f] = h$ is straightforward. Let us now check that $1$
is the highest weight vector:
$$ h 1 = (td) 1 - (t^{-1} d) 1 = (\alpha-x) - (-\alpha-x) = (2 \alpha) 1 ,$$
$$ e 1 = {1\over 2}  \left( (td)1  - 2d 1 + (t^{-1} d)1 \right) 
= {1 \over 2}  \left( (\alpha-x) + 2x + (-\alpha-x) \right) = 0 .$$
 Finally, we have noted above that $\C_\alpha[x]$ is generated by its highest weight vector $1$. Thus it is 
isomorphic to the Verma module over $\sll_2 (\C)$.

\

{\bf Corollary 4.2.} 
{\it
The action of $\sll_2(\C)$ on its Verma module $M_\lambda$ can be extended to the action of the
Lie algebra of vector fields on a circle, which is irreducible except when $\lambda = 0$.
}

{\bf Proof.} The Verma module $M_\lambda$ over $\sll_2$ can be realized as $\C_\alpha [x]$ with $\alpha = {\lambda \over 2}$.
It follows from (\Wact) that the action of $\sll_2$ may be extended to the Lie algebra $\Der \C[t,t^{-1}]$. Let us show
that the module $\C_\alpha [x]$ is irreducible over the Lie algebra of vector fields on a circle, except when $\alpha = 0$.
It is well-known that when $\lambda$ is a non-negative integer, the Verma module over $\sll_2$ has a single proper submodule 
of codimension $\lambda +1$, and 
for all other $\lambda$ the Verma module is irreducible. To prove that $\C_\alpha [x]$
is irreducible as a module over $\Der \C[t,t^{-1}]$, all we need to show is that the proper $\sll_2$-submodule is not invariant     
under the action of the Lie algebra of vector fields on a circle. Assume the contrary. Then we would get a finite-dimensional
quotient module. But the Lie algebra of vector fields on a circle, being a simple infinite-dimensional algebra, can't have 
non-trivial finite-dimensional representations. Thus $\C_\alpha [x]$ is irreducible as a module over this Lie algebra,
unless $\alpha = 0$. In the case when $\alpha = 0$, the action (\Wact) on $\C_0[x]$ is
$$(t^m d) p(x) = -x p(x+m) $$
and we can immediately see that the space of polynomials vanishing at $0$ forms a 
submodule in $\C_0[x]$.


\

 Let us now consider the action of the Casimir operator
$$ \Omega = - {1\over 2} (t d) \otimes (t^{-1} d) +  d \otimes d - {1\over 2} (td) \otimes (t^{-1} d) 
= {1\over 2} e \otimes f + {1\over 4} h \otimes h + {1\over 2} f \otimes e $$
on the tensor product of modules $\C_{\alpha_1} [x] \otimes \C_{\alpha_2} [y] \cong M_{2\alpha_1} \otimes M_{2\alpha_2}$.
It is well-known that $\Omega$ commutes with the action of $\sll_2(\C)$ on the tensor product of these modules.

 Let us write down the action of $\Omega$ on 
$\C_{\alpha_1} [x] \otimes \C_{\alpha_2} [y]$ explicitly:
$$\Omega p(x,y) = - {1\over 2} (x - \alpha_1) (y + \alpha_2) p(x+1,y-1)
+ xy p(x,y) - {1\over 2} (x + \alpha_1) (y - \alpha_2) p(x-1,y+1) .$$
Although each term in the right hand side increases the degree of $p(x,y)$,
it can be seen that the Casimir operator $\Omega$ actually preserves the space
of polynomials of total degree in $x$ and $y$ at most $n$.

 The Casimir operator $\Omega$ plays an important role in the differential equations that we study in the next section.
In particular we would need to know the spectrum of $\Omega$ on $\C_{\alpha_1} [x] \otimes \C_{\alpha_2} [y]$.

\

{\bf Proposition 4.3.} 
{\it
The space $\C_{\alpha_1} [x] \otimes \C_{\alpha_2} [y]$ decomposes into a direct sum of 
finite-dimensional $\Omega$-invariant subspaces. The eigenvalues of $\Omega$ are 
$$\left\{ \omega_n = {n^2 \over 2} - n(\alpha_1 + \alpha_2 + {1\over 2}) + \alpha_1 \alpha_2 \big| n=0,1,2,\ldots \right\}$$ 
(each occurring with infinite multiplicity).
}

{\bf Proof.} Since the Casimir operator has weight $0$, it preserves the subspaces
$V_n = \mathop\oplus\limits_{k = 0}^n \C f^k 1_{\alpha_1} \otimes f^{n-k} 1_{\alpha_2}$. This establishes 
the first claim of the proposition. Next, the operator $f \otimes \id + \id \otimes f$ is an injection
of $V_{n-1}$ into $V_n$. Since $\Omega$ commutes with this operator (it is the action of $f$ on the tensor
product of two modules), we conclude that $V_n$ contains an $\Omega$-invariant subspace of codimension $1$,
on which the spectrum of $\Omega$ coincides with its spectrum on $V_{n-1}$.

Hence there exists a sequence $\left\{ \omega_n | n = 0,1,2, \ldots \right\}$ such that the spectrum of $\Omega$ on 
$V_n$ is $\left\{ \omega_0, \ldots, \omega_n \right\}$. The sum $\omega_0 + \ldots + \omega_n$ equals the trace of
$\Omega$ on $V_n$, which is easy to compute, since on these subspaces $\tr (\Omega) = {1\over 4} \tr (h \otimes h)$
(the other two terms are traceless). We have
$$ \omega_0 + \ldots + \omega_n = {1\over 4} \sum_{k=0}^n (-2k+ 2\alpha_1)(-2n+2k+2\alpha_2), $$
and thus for $n\geq 1$,
$$\omega_n = \sum_{k=0}^n (k-\alpha_1)(n-k-\alpha_2) -  \sum_{k=0}^{n-1} (k-\alpha_1)(n-1-k-\alpha_2)
= -(n-\alpha_1) \alpha_2 +  \sum_{k=0}^{n-1} (k-\alpha_1) $$
$$= {n^2 \over 2} - n(\alpha_1 + \alpha_2 + {1\over 2}) + \alpha_1 \alpha_2 . \eqno{(\spec)}$$
Note that this formula also gives the correct value of $\omega_0 = \alpha_1 \alpha_2$.

{\bf Remark.} If all $\omega_n$ are distinct (which happens when the sum of the highest weights
$2(\alpha_1 + \alpha_2)$ is not a non-negative integer),
then $\Omega$ is diagonalizable on the tensor product of these Verma modules, and the tensor product decomposes into a direct
sum of Verma modules over $\sll_2(\C)$.

\

\

{\bf 5. Differential equations in modules for vertex algebras.}

\

 Our goal is to show that the modules $L_\H (\alpha, \gamma)$ with $\alpha\not\in\Q$,
are generated by $T$ as the modules over the loop subalgebra $\L$. This will give a complete 
description of these modules. We will do this by studying a differential equation in 
$L_\H (\alpha, \gamma)$ that arises from (\Ydom).

 The first step towards this result is the following

{\bf Lemma 5.1.} 
{\it
Let us suppose that for all $j$, 
$(t_0^{-j} t_1^m d_1) v(n)$ belongs to the subspace $U(\L_-) T$ in the module $L_\H (\alpha, \gamma)$.
Then $L_\H (\alpha, \gamma) =  U(\L_-) T$.
}

{\bf Proof.} The proof is based on a simple Poincar\'e-Birkhoff-Witt argument.
The module $L_\H (\alpha, \gamma)$ is spanned by the vectors 
$$ (t_0^{-j_1} t_1^{m_1} d_1) \ldots (t_0^{-j_k} t_1^{m_k} d_1) v(n) , \quad j_1, \ldots, j_k > 0.$$
We are going to show that these vectors belong to the space $U(\L_-) T$ using a double
induction on total degree $-j=-(j_1+\ldots +j_k)$ and on length $k$.
If the total degree is zero, then there is nothing to prove. Suppose now that the claim of the lemma
is true for all degrees above $-j$. We have that
 $ (t_0^{-j_2} t_1^{m_2} d_1) \ldots (t_0^{-j_k} t_1^{m_k} d_1) v(n)$ has degree greater than $-j$.
By induction assumption, it belongs to the space $U(\L_-) T$. Thus without the loss of generality
we may assume that $|m_2|, \ldots, |m_k| \leq 1$.

Let us prove by induction on length $k$ that for the claim of the
lemma holds for the vectors in $L_\H (\alpha, \gamma)$ of degree $-j$. The basis of induction, $k=1$, is
precisely the assumption of the lemma. Let $k \geq 2$. Then
$$ (t_0^{-j_1} t_1^{m_1} d_1) (t_0^{-j_2} t_1^{m_2} d_1) \ldots (t_0^{-j_k} t_1^{m_k} d_1) v(n)$$
$$= (m_2 - m_1) (t_0^{-j_1-j_2} t_1^{m_1+ m_2} d_1) \ldots (t_0^{-j_k} t_1^{m_k} d_1) v(n)$$
$$+ (t_0^{-j_2} t_1^{m_2} d_1) (t_0^{-j_1} t_1^{m_1} d_1) \ldots (t_0^{-j_k} t_1^{m_k} d_1) v(n).$$
The first summand in the right hand side belongs to the space $U(\L_-) T$ by induction assumption 
on length, while for the second summand we get that 
$$ (t_0^{-j_1} t_1^{m_1} d_1) (t_0^{-j_3} t_1^{m_3} d_1) \ldots (t_0^{-j_k} t_1^{m_k} d_1) v(n) \in U(\L_-) T$$
by the induction assumption on degree. Noting that $t_0^{-j_2} t_1^{m_2} d_1 \in \L_-$
since $|m_2|\leq 1$, we obtain the desired claim. The lemma is now proved.

\

 Thus it remains to prove that $(t_0^{-j} t_1^m d_1) v(n) \in U(\L_-) T$ for all $j$. We will achieve 
this by studying the generating series 
$$ P(m,n,z) = d_1(m,z) v(n) = z^{-1} (t_1^m d_1) v(n) + (t_0^{-1} t_1^m d_1) v(n) + z (t_0^{-2} t_1^m d_1) v(n) + \ldots  $$

 Since $d_1(m,z)$ satisfies the differential equation (\Ydom), we conclude that $P(m,n,z)$ satisfies
an analogous differential equation
$$\displaylines{ {\d \over \d z} P(m,n,z)  =  \cr
 \left( - {m+1\over 2} \lcol d_1(1,z) d_1(m-1,z)\rcol   
+ m \lcol d_1(0,z) d_1(m,z)\rcol \right. \hfill \cr  
\hfill \left. - {m-1\over 2} \lcol d_1(-1,z) d_1(m+1,z)\rcol \right) v(n) \cr
 = \left( - {m+1\over 2} d_1(1,z)_+ d_1(m-1,z)   
+ m d_1(0,z)_+ d_1(m,z) \right. \hfill\cr 
\hfill \left. - {m-1\over 2} d_1(-1,z)_+ d_1(m+1,z) \right)  v(n) \cr
+z^{-1} \left( - {m+1\over 2} d_1(m-1,z) (t_1 d_1) 
+ m z^{-1} d_1(m,z) d_1 \right. \hfill\cr
\hfill \left. - {m-1\over 2} z^{-1} d_1(m+1,z) (t_1^{-1}d_1) \right) v(n) \cr
 = - {m+1\over 2} d_1(1,z)_+ P(m-1,n,z)  
+ m d_1(0,z)_+ P(m,n,z) 
- {m-1\over 2} d_1(-1,z)_+ P(m+1,n,z) \cr
\hfill +z^{-1} \Omega P(m,n,z) , \hfill \quad (\vdom) \cr}  $$
where
$$d_1(k,z)_+ = \sum_{j<0} (t_0^j t_1^k d_1) z^{-j-1} \eqno{(\nopm)}$$
and
$$\displaylines{
 \Omega P(m,n,z) = \cr
\hfill - {1\over 2} (m+1)(n-\alpha) P(m-1,n+1,z)
+   mn P(m,n,z) - {1\over 2} (m-1)(n+\alpha) P(m+1,n-1,z) .  
\hfill {(\vom)} \cr} $$

In order to solve this differential equation, we will treat $m$ and $n$ as formal variables.
From this point of view, the transformation $\Omega$ defined in (\vom), 
is essentially the
Casimir operator on the space $\C_{\alpha_1}[m] \otimes \C_{\alpha_2} [n]$ with $\alpha_1 = 1$, $\alpha_2 = \alpha$.  

 Let us give a precise description of the spaces which we will be working with. Consider the algebra $\C[m,n]$ of polynomials in two
formal variables $m,n$.
Let $\Gamma$ be the space of affine functions
in $m$ with integer coefficients: $\Gamma=\{ am + b | a,b \in\Z \}$. We are going to replace in the definitions of
Lie algebras and modules which we have previously discussed, the algebra of Laurent polynomials $\C[t_0^{\pm 1}, t_1^{\pm 1}]$,
which is a group algebra of $\Z \oplus \Z$, with the group algebra of $\Z \oplus \Gamma$, and extend the scalars to $\C[m,n]$: 
$$\hR = \C[m,n] \otimes \C [\Z \oplus \Gamma],$$ 
$$ \hH = \hR d_1, $$
with obvious modifications of the formulas (\Lieb) for the Lie bracket:
$$ [t_0^i t_1^\mu d_1, t_0^j t_1^\nu d_1] = (\nu - \mu) t_0^{i+j} t_1^{\mu+\nu} d_1, \quad i,j \in \Z, \mu,\nu\in \Gamma.$$
In the same way we will replace the $\H_0$-module $T(\alpha,\gamma)$ with the $\hH_0$-module $\hT(\alpha,\gamma)$,
$\alpha,\gamma \in\C$:
$$\hT = \mathop\bigoplus_{\nu \in \Gamma} \C[m,n] v(n+\nu) .$$
As before, we form the generalized Verma modules $M_\hH (\hT) = U(\hH_-) \otimes_{\C[m,n]} \hT$ and $M_\H (\hT) = U(\H_-) \otimes_{\C} \hT$. 
%
%
%
We define $L_\hH (\hT)$ (resp. $L_\H (\hT)$) to be the quotient of
$M_\hH (\hT)$ (resp. $M_\H (\hT)$) 
by the maximal $\Z \oplus \Gamma$-graded submodule trivially intersecting 
$\hT$.

 The following lemma is quite obvious:

{\bf Lemma 5.2.} 
{\it
(a) A specialization  $m \mapsto m_0 \in \Z$, $n \mapsto n_0 \in \Z$
extends to a homomorphism of Lie algebras $\hH \rightarrow \H$ and a homomorphism of modules 
$$M_\hH (\hT) \rightarrow M_\H (T).$$

(b) The above map factors through to the following:
$$L_\hH (\hT) \rightarrow L_\H (T).$$
}

 Our strategy will be to solve the differential equation (\vdom) in the module $L_\hH (\hT)$, where we treat $m$ and $n$
as formal variables, and then use part (b) of the Lemma above to obtain a solution in $L_\H (T)$.

 We also need to establish a relation between $M_\H (\hT)$ and $L_\hH (\hT)$:

{\bf Lemma 5.3.} 
{\it
The composition 
$$ \varphi: \quad M_\H (\hT) \rightarrow L_\hH (\hT) $$
of an embedding $M_\H (\hT) \rightarrow M_\hH (\hT)$ with a projection
$M_\hH (\hT) \rightarrow L_\hH (\hT)$ is surjective.
}

{\bf Proof.} We need to show that for every $x \in M_\hH (\hT)$ there exists a $y \in M_\H (\hT)$, such that 
the images of $x$ and $y$ in $L_\hH (\hT)$ coincide. It is sufficient to consider the case when $x$ is a 
monomial:
$$ x = (t_0^{-j_1} t_1^{\mu_1} d_1) \ldots (t_0^{-j_s} t_1^{\mu_s} d_s) v(n+\nu), \quad \mu_1,\ldots,\mu_s, \nu \in \Gamma,
j_1, \ldots j_s > 0.$$
The kernel of the projection $M_\hH (\hT) \rightarrow L_\hH (\hT)$ is invariant under the shifts
$n \mapsto n + \mu, \, \mu\in \Gamma$. Thus we may assume without the loss of generality that 
$\mu_1+\ldots+\mu_s + \nu = 0$. A homogeneous element of $M_\hH (\hT)$ is non-zero in $L_\hH (\hT)$ 
if and only if a submodule generated by this element in $M_\hH (\hT)$ intersects with $\hT$ non-trivially. Thus the image
of $x$ in $L_\hH (\hT)$ is determined by the values of $u x$, where $u$ are elements in $U(\hH_+)$ of degree $k = j_1 + \ldots + j_s$
in $t_0$. Moreover, it is sufficient to take $u$ to be elements of degree $k$ in $U(\H_+)$. This follows from the fact that a 
non-zero polynomial in $m$ assumes a non-zero value at some integer. 
We may also assume that $u$ is also homogeneous in $t_1$ of degree $i$.
Then
$$ u x = f (m,n) v(n+i),$$
where $f$ is a polynomial of degree at most $s+k$ in $m$. However, a polynomial of degree $s+k$ may be extrapolated through
its $s+k+1$ distinct integer values. Hence we can construct $y$ to be an extrapolation in $m$ through the specializations
of $x$ to any $s+k+1$ distinct integer values of $m$. Note that such an extrapolation of $x$ depends only on the degree
$s+k$, but not on a particular raising element $u$.

\

{\bf Example.} The same technique may be applied to the vertex algebras $V_\H$ and $\bV_\H$. Consider 
the element $x = (t_0^{-1} t_1^m d_1)(t_0^{-1} t_1^n d_1) \o \in V_\H$. We would like to get a reduction rule for
the image of $x$ in $\bV_\H$. First, we make a shift $n \mapsto n-m$, and get 
$x^\prime = (t_0^{-1} t_1^m d_1)(t_0^{-1} t_1^{n-m} d_1) \o$.
The image of $x^\prime$ in $\bV_\H$ is determined by the values of the raising operators of degree one, 
$u = t_0 t_1^i d_1, \, i\in\Z$, applied to $x^\prime$. It can be easiy seen that the coefficients that occur is $ux^\prime$
are polynomials in $m$ of degree two (cf., last formula in (\rais)). Thus $x^\prime$ may be interpolated in $\bV_\H$
with a quadratic polynomial in $m$ using any three integer specializations of $m$ in $x^\prime$.
Choosing $-1,0,1$ as the interpolation points, we get (cf., (\drm)):
$$(t_0^{-1} t_1^m d_1)(t_0^{-1} t_1^{n-m} d_1) \o = 
{m(m+1) \over 2} (t_0^{-1} t_1 d_1)(t_0^{-1} t_1^{n-1} d_1) \o $$
$$-(m-1)(m+1) (t_0^{-1} d_1)(t_0^{-1} t_1^{n} d_1) \o
+ {m(m-1) \over 2} (t_0^{-1} t_1^{-1} d_1)(t_0^{-1} t_1^{n+1} d_1) \o .$$
Taking the shift back, $n \mapsto n+m$, we obtain analogous relation for $x$.

\

Let us consider a subspace $W$ in $L_\hH (\hT)$, spanned as a $\C[m,n]$-module by the elements
$(t_0^{-j_1} t_1^{\mu_1} d_1) \ldots (t_0^{-j_s} t_1^{\mu_s} d_s) v(n+\nu)$, \quad  with
$\mu_1+\ldots+\mu_s + \nu = m + i, \quad i\in\Z$, and the corresponding subspace $\tW$ in $M_\H (\hT)$:
$$\tW = U(\H_{-}) \otimes \left( \mathop\bigoplus_{i\in\ZS}  \C[m,n] v(n+m+i) 
\right) .$$
We note that $\varphi(\tW) = W.$

 The coefficients at all powers of $z$ in $P(m,n,z)$ belong to the space $W$, and the same is true for all summands
in the right hand side of the differential equation (\vdom).
This allows us to restrict to subspace $W$ when solving (\vdom).
The reason for considering such a restriction is that in space $\tW$ the Casimir operator $\Omega$ acts just on the coefficients:
$$\Omega_\tW \left( f(m,n) u v(m+n+i) \right) =  \Omega \left( f(m,n) \right) u v(m+n+i),\quad u \in U(\H_{-}), i \in\Z. \eqno{(\otw)} $$
Since the kernel of $\varphi$ is invariant under the shift operators $m \mapsto m+j, n \mapsto n+i, \,
i,j \in \Z$, it is also invariant under $\Omega$. Thus we may calculate the action of $\Omega$ in $W$ by
lifting elements to $\tW$: 
$$\Omega_W = \varphi \Omega_\tW \varphi^{-1} . \eqno{(\opw)} $$

 As a consequence, we get that the action of $\Omega_W$ is locally finite, i.e., every vector in $W$
is contained in a finite-dimensional  $\Omega$-invariant $\C$-subspace, and that the spectrum of 
$\Omega_W$ is a subset of (\spec) with $\alpha_1 = 1, \alpha_2 = \alpha$.
 
Denote by $S_i$ the shift operator $m \mapsto m+i$, $i \in \Z$, on the space $W$. Then the differential equation (\vdom) may be rewritten as
follows:
$${\d \over \d z} P(m,n,z) = \left( z^{-1} \Omega + A_+(z) \right) 
P(m,n,z), \eqno{(\szp)} $$
where $\Omega$ is the Casimir operator (\vom) and
$$A_+(z) =  - {m+1\over 2} d_1(1,z)_+ S_{-1}  
+ m d_1(0,z)_+ - {m-1\over 2} d_1(-1,z)_+ S_1 . $$

\

\

{\bf 6. Solving non-commutative differential equations.}

\

 Consider a differential equation in a vector space $W$:
$$ {\d \over \d z} w(z) = A(z) w(z),  \eqno{(\diffeq)}$$
where 
$$ w(z) = \sum_{j \geq r} w_{j} z^j, \quad \hbox{\rm with \ } w_{j} \in W, \eqno{(\eqw)}$$
and
$$ A(z) = \Omega z^{-1} + \sum_{i \geq 0} A_{i} z^i, \quad \hbox{\rm with \ } 
\Omega, A_{i} \in \End(W). \eqno{(\eqA)}$$
We would like to solve this differential equation for $w(z)$ given the ''initial data''
$w_{r} \in W$. We will denote by $A_+(z)$ the tail of $A(z)$:
$$A_+(z) = \sum_{i \geq 0} A_{i} z^i .$$

In order to solve this differential equation, we will introduce a special type of integration.
An expression $b(z) = \sum\limits_{j=-\infty}^\infty b_{j} z^j$, involving $z^{-1}$, can't be integrated in a usual way, 
since $z^{-1}$ does not have  an antiderivative in power series. To avoid this difficulty we define
$$ z^\Omega \int z^{-\Omega} b(z) dz = \sum_{j=-\infty}^\infty \left( (j+1)I - \Omega \right)^{-1} b_{j} z^{j+1} . \eqno{(\zint)}$$
For this integral to be defined we require that the operator $(j+1)I - \Omega$ is invertible whenever $b_{j} \neq 0$.

 The notation for this integral is motivated by the formula
$$ z^\omega \int z^{-\omega} b(z) dz = \sum_{j=-\infty}^\infty {1 \over j+1 - \omega} b_{j} z^{j+1} ,$$
for $\omega \in\C, \omega\not\in\Z$.

{\bf Lemma 6.1.} 
{\it
Assume that the integral (\zint) is well-defined. Then
$${\d \over \d z} z^\Omega \int z^{-\Omega} b(z) dz
= z^{-1} \Omega \left( z^\Omega \int z^{-\Omega} b(z) dz \right) + b(z) .$$
}

{\bf Proof.} 
$$  {\d \over \d z} z^\Omega \int z^{-\Omega} b(z) dz $$
$$ = \sum_{j=-\infty}^\infty (j+1) \left( (j+1)I - \Omega \right)^{-1} b_{j} z^{j} $$
$$ = \sum_{j=-\infty}^\infty \left( (j+1)I - \Omega \right) \left( (j+1)I - \Omega \right)^{-1} b_{j} z^{j} $$
$$ + \sum_{j=-\infty}^\infty \Omega \left( (j+1)I - \Omega \right)^{-1} b_{j} z^{j} $$ 
$$ = z^{-1} \Omega \left( z^\Omega \int z^{-\Omega} b(z) dz \right) + b(z) .$$

 The following theorem gives a formula for the solution of the differential equation (\diffeq).

{\bf Theorem 6.2.} 
{\it
Fix $r\in \Z$.
Consider a non-commutative differential equation 
$$ {\d \over \d z} w(z) = A(z) w(z), $$
with $w(z)$ and $A(z)$ as in (\eqw) and (\eqA).  
Assume that the operators $\Omega - jI$ are invertible for all $j > r$. Given $w_{r} \in W$ satisfying the consistency condition
$$ \left( \Omega - rI \right) w_{r} = 0 , \eqno{(\cons)}$$
there exists a unique solution of the above differential equation of the form
$$ w(z) = w_{r} z^r + \sum\limits_{j > r} w_{j} z^j, \quad w_{j} \in W.$$
This solution is given by the formula
$$ w(z) = \sum_{k=0}^\infty P_k (z) \eqno{(\Pnz)}$$
with $P_0(z) = w_{r} z^r$ and
$$P_{k+1} (z) = z^\Omega \int z^{-\Omega} A_+(z) P_k(z) dz 
\quad \hbox{\rm for all \ } k\geq 0. \eqno{(\Poz)}$$
The series $P_k(z)$ has no powers of $z$ below $r+k$, which makes the infinite sum (\Pnz) well-defined.
}

{\bf Proof.} Let us expand the differential equation (\diffeq) in powers of $z$. We can see that it does not have
any terms in powers of $z$ below $r-1$. The $z^{r-1}$ term of this equation is 
$$ r w_{r} = \Omega w_{r}, $$
which holds by (\cons). The term at $z^{j-1}$ with $j > r$ is
$$ j w_{j} = \Omega w_{j} + \sum_{i=r}^{j-1} A_{j-i-1} w_{i}.$$
By the assumption of the theorem, the operators $\Omega - jI$ with $j > r$ are invertible,
which gives the recursive formula for the (unique) solution of (\diffeq):
$$ w_{j} = \left( jI - \Omega \right)^{-1} \sum_{i=r}^{j-1} A_{j-i-1} w_{i}. \eqno{(\wj)}$$

 Finally, let us prove that the same solution is given by the formula (\Pnz).
To do this, we need to show that the series defined by (\Pnz), (\Poz) also satisfies 
(\diffeq). 

Since $A_+(z)$ contains only non-negative powers of $z$, we see that integration (\zint) increases the powers of $z$ in (\Poz) by 1. Thus $P_k(z)$ has no powers of $z$ below $r+k$, which makes the infinite sum (\Pnz) well-defined.

By Lemma 6.1, 
$$ {\partial \over \partial z} P_{k+1} (z) = z^{-1} \Omega P_{k+1} (z) + A_+(z) P_k(z). $$
Also by (\cons) we get that 
$$ {\partial \over \partial z} P_{0} (z) = z^{-1} \Omega P_{0} (z).$$
Thus
$$  {\partial \over \partial z} \sum_{k=0}^\infty P_k (z) = 
z^{-1} \Omega \sum_{k=0}^\infty P_k (z) + A_+(z) \sum_{k=1}^\infty P_{k-1} (z) ,$$
and we see that $\sum_{k=0}^\infty P_k (z)$ satisfies (\diffeq). Since the solution of
(\diffeq) is unique, we conclude that the solutions (\wj) and (\Pnz)-(\Poz) coincide.

\

 We are now going to establish our main result, Theorem 1.2, by solving the differential equation (\szp) in the space $W \subset L_\hH (\hT)$. Let us first verify that (\szp) satisfies the assumption of Theorem 6.2. By (\otw), (\opw) and Proposition 4.3, the
action of the Casimir operator $\Omega$ on $W$ is locally finite, and it has the spectrum
$$ \left\{ \omega_n = {1\over 2} \left( n^2 - (2\alpha + 3)n + 2\alpha \right) | n=0,1,2, \ldots \right\} .$$
It can be easily seen that $\omega_1 = -1$ for all $\alpha$, whereas
the spectrum of $\Omega$ contains a non-negative integer if and only if 
$\alpha \in \Q$. Hence for all $\alpha \not\in\Q$ the operators 
$\Omega - jI$ are invertible for all $j = 0,1,2,\ldots$ 

 The first term in the series $P(m,n,z)$ is 
$$P_0 (m,n,z) = z^{-1} (t^m d) v(n) = z^{-1} (n - \alpha m) v(n+m).$$
Thus the value of $r$ that appears in the statement of Theorem 6.2 is
$r = -1$. One can verify that the consistency condition
$$(\Omega + I) (n - \alpha m) v(n+m) = 0 $$
holds -- and actually this follows from the fact that (\vdom) does have a solution in $W$.

 Hence, by Theorem 6.2, the solution of (\vdom) is given by the formula
$$P(m,n,z) = \sum_{k=0}^\infty P_k (m,n,z),$$
with
$$\displaylines{
P_{k+1} (m,n,z) = z^\Omega \int z^{-\Omega} A_+(z) P_k(m,n,z) dz \cr
= z^\Omega \int z^{-\Omega} 
\left(
-{m+1 \over 2} d_1 (1,z)_+ P_k(m-1,n,z) + m d_1 (0,z)_+ P_k(m,n,z) \right.
\hfill \cr
\hfill \left. -{m-1 \over 2} d_1 (-1,z)_+ P_k(m+1,n,z) \right) dz .\cr}$$

 Since $P_0(m,n,z)$ belongs to $\hT$ and the moments of 
$d_1 (1,z)_+, d_1 (0,z)_+, d_1 (-1,z)_+$ are in the loop algebra $\L$,
we conclude that the moments of $P(m,n,z) = d_1 (m,z) v(n)$, viewed as
elements of $L_\hH (\hT)$, belong in fact to the subspace $U(\L_-) \hT$.
Applying Lemma 5.2 (b), we obtain that for all $m,n,j \in \Z$,
$$(t_0^j t_1^m d_1) v(n) \in U(\L_-) T$$
in the module $L_\H (T)$. Thus by Lemma 5.1, we get that $L_\H (T)$
is generated by $T$ as an $\L$-module,
$$L_\H (T) = U(\L_-)T.$$
Hence we have a surjective map of $\L$-modules
$$M_\L (T) \rightarrow L_\H (T).$$

 In the Appendix we use results of [F] and [KM] to show that the generalized Verma module 
$M_\L (\alpha, \gamma)$ over the loop algebra 
$\L$ has no non-trivial submodules trivially intersecting with $T$
whenever $\alpha \not\in {1\over 2} \Z$. Since 
$\alpha \not\in \{ -1, 0 \}$, $T$ is an irreducible $\H_0$-module, and thus for 
$\alpha \not\in \Q$ the kernel of the map 
$M_\L (T) \rightarrow L_\H (T)$ is trivial, and we obtain an 
%
%
isomorphism of $\L$-modules:
$$L_\H (T) \cong M_\L(T) .$$
Finally, by Theorem 3.8, for $\beta = - {\alpha(\alpha+1) \over 2}$, the modules 
$L_\D(\alpha,\beta,\gamma)$ and $L_\H (\alpha, \gamma)$ are isomorphic as $\H$-modules, 
and hence as $\L$-modules. This gives us 
$$L_\D(\alpha,\beta,\gamma) \cong L_\H (\alpha, \gamma) \cong 
M_\L (\alpha, \gamma) .$$
We conclude that for $\alpha \not\in\Q$, $\beta  = - {\alpha(\alpha+1) \over 2}$, the module 
$L_\D(\alpha,\beta,\gamma)$ remains irreducible as a module over the loop subalgebra $\L$, 
and its action on 
$M_\L (\alpha, \gamma)$ can be extended to the larger algebra $\D$ of vector fields 
on a 2-dimensional torus.

\

\

{\bf Appendix. Generalized Verma modules for the loop Lie algebra.}

\

 Futorny [F] defined the Shapovalov form on the generalized Verma module
$M_\L (\alpha, \gamma)$, and used it to establish a criterion for irreducibility of $M_\L (\alpha, \gamma)$ (in fact [F] treats a more general case of modules for the affine Kac-Moody algebra $\widehat {sl}_2$, however we restrict our attention here to the level zero case of the loop algebra $\L$). A complete formula for the Shapovalov determinant was obtained by Khomenko and Mazorchuk in [KM]. Here we would like to address a related question of existence of non-trivial submodules in 
$M_\L (\alpha, \gamma)$ having a trivial intersection with the top
$T (\alpha, \gamma)$. We will do this with the help of the Shapovalov determinant formula given in [KM].

 Consider three related partition functions,
$$\eqalign{
p_+ (s,n) &= \dim U (\L_+)_{(s,n)} , \cr
p (s,n) &= \dim U (\C (t_1 d_1) \oplus \L_+)_{(s,n)} , \cr
p_3 (s) &= \dim U (\L_+)_s . \cr
}$$
Here a $\Z^2$-grading refers to the grading by degrees in $t_0$ and $t_1$, while a $\Z$-grading is taken just by degrees in $t_0$.

 From the structure of the root system of $\L$ we see
that $p_+(s,n) = 0$ when $|n| > s$ and we have the symmetry
$p_+(s,n) = p_+ (s, -n)$. The Poincar\'e-Birkhoff-Witt theorem implies
that 
$$p(s,n) = \sum_{k\geq 0} p_+ (s, n-k), \eqno{(\PPsn)} $$
while 
$$p_3(s) = \sum_{n\in\ZS} p_+ (s, n). \eqno{(\PTsn)}$$
We also note that the partition function $p_3$ appears in the expansion of
(\charf):
$$  \prod\limits_{n\geq 1} (1-q^{n})^{-3} = \sum_{s=0}^\infty p_3 (s) q^s .$$

 Let us define two determinants associated with the module $M_\L (\alpha, \gamma)$. As before, we consider $\Z^2$-grading of $M_\L (\alpha, \gamma)$ by degrees in $t_0$ and $t_1$. We will assume that the top $T(\alpha,\gamma)$ has degree zero in $t_0$. Note that any submodule in $M_\L (\alpha, \gamma)$ is homogeneous with respect to $\Z^2$-grading since
$d_1 \in \L$, while the action of $d_0$ can be recovered via the Sugawara construction (cf., (\Ydzm)).

Fix homogeneous bases $\{ u_i^+ | i = 1, \ldots p_3(s) \}$ in
$U(\L_+)_s$ and $\{ u_i^- \}$ in $U(\L_-)_{-s}$. Let $r_i$ be the degree
in $t_1$ of $u_i^+$. We may assume that the degree in $t_1$ of $u_i^-$ is
$-r_i$. Define $\epsilon_i = 1$ if $r_i >0$ and $\epsilon_i = -1$ if $r_i < 0$. The value of $\epsilon_i$ when $r_i = 0$ will not be relevant.

In this set-up, the set $\{ u_i^- v(n+r_i) \}$ forms a basis of
$M_\L (\alpha, \gamma)_{(-s,n)}$, where $s \geq 1$, $n \in \gamma + \Z$.

 We introduce two square matrices. Define $G_{ij} = G_{ij}^{s,n}$ 
to be the coefficient at 
\break
$v(n+r_i)$ in $u_i^+ u_j^- v(n+r_j)$.
Similarly, define $F_{ij} = F_{ij}^{s,n}$ to be the coefficient at $v(n)$
in $(t_1^{-\epsilon_i} d_1)^{|r_i|} u_i^+ u_j^- (t_1^{\epsilon_j} d_1)^{|r_j|} v(n)$. 
We will consider the determinants of these matrices as polynomials in $n$ and $\alpha$.

 It is clear that $M_\L (\alpha, \gamma)$ has a non-trivial submodule that trivially intersects the top $T(\alpha,\gamma)$ if and only if 
$\det G^{s,n} = 0$ for some $s \geq 1, n \in \gamma+\Z$. The second determinant, $\det F$, is the determinant of the Shapovalov form, calculated in [KM]. Since there is a simple connection between the two determinants, we will get a formula for $\det G$.

 Indeed, since
$$(t_1^{\epsilon_j} d_1)^{|r_j|} v(n) = 
\prod_{k=1}^{|r_j|} (n+(k-1)\epsilon_j - \alpha \epsilon_j) v(n+r_j),$$
and similarly
$$(t_1^{-\epsilon_i} d_1)^{|r_i|} v(n+r_i) = 
\prod_{k=1}^{|r_i|} (n+k\epsilon_i + \alpha \epsilon_i) v(n),$$
we get that
$$\det F = \det G \times 
\prod_{j=1}^{p_3(s)} \prod_{k=1}^{|r_j|} (n+(k-1)\epsilon_j - \alpha \epsilon_j)
 \times  \prod_{i=1}^{p_3(s)} 
 \prod_{k=1}^{|r_i|} (n+k\epsilon_i + \alpha \epsilon_i) . 
\eqno{(\FGrel)}$$
Combining the factors in (\FGrel) that correspond to the basis elements of equal positive or negative degrees, we obtain
$$\det F = \det G \times
\prod_{m=1}^s \prod_{k=1}^m \big(
(n+k-1-\alpha)(n-k+1+\alpha)(n+k+\alpha)(n-k-\alpha) \big)^{p_+ (s,-m)}.$$
Finally, interchanging the order of the products and taking 
(\PPsn) into account, we get
$$\det F = \det G \times
 \prod_{k=1}^s \big(
(n+k-1-\alpha)(n-k+1+\alpha)(n+k+\alpha)(n-k-\alpha) \big)^{p (s,-k)}.$$

The Shapovalov determinant formula given in [KM], applied to our situation yields (up to a non-zero constant factor):
$$
\displaylines{
\det F = 
\prod_{k=1}^s \big(
(n+k-1-\alpha)(n-k+1+\alpha)(n+k+\alpha)(n-k-\alpha) \big)^{p (s,-k)} \cr
\times
\prod_{m=1}^s \prod_{k=1}^{\left[ {s \over m} \right]}
\big( (2m-k+2\alpha+1) (2m-k-2\alpha-1) \big)^{p_3 (s-mk)} .
\cr} $$ 
 
Comparing the last two equalities, we obtain a formula for $\det G$
(up to a non-zero constant factor):
$$\det G^{s,n} = 
\prod_{m=1}^s \prod_{k=1}^{\left[ {s \over m} \right]}
\big( (2m-k+2\alpha+1) (2m-k-2\alpha-1) \big)^{p_3 (s-mk)} .$$

As an immediate corollary, we get

{\bf Proposition A.1.} 
{\it
The $\L$-module $M_\L (\alpha, \gamma)$ has a non-trivial submodule trivially intersecting the top $T(\alpha, \gamma)$ if and only if $\alpha \in {1\over 2} \Z$.
}

\

{\bf Remark.} The above argument can be applied in the general case considered in [KM]. For the generalized Verma modules 
$M_{\widehat \L} (\alpha, \gamma, c)$ for affine Kac-Moody algebra
${\widehat \L} = {\widehat {sl}}_2$ at level $c$, we get

$$\det G^{s,n} = 
\prod_{m=1}^s \prod_{k=1}^{\left[ {s \over m} \right]}
\big( (c+2)(m(c+2)-k+2\alpha+1) (m(c+2)-k-2\alpha-1) \big)^{p_3 (s-mk)} .$$
As a consequence, we get that $M_{\widehat \L} (\alpha, \gamma, c)$ 
has a non-trivial submodule trivially intersecting the top $T(\alpha, \gamma)$ if and only if $2 \alpha +1 = \pm (m(c+2) -k)$ for some integer
$m,k \geq 1$ or $c=-2$ (cf., [F], Theorem 3.11).

\

\

{\bf References:}

\

\noindent
[BB] Berman, S., Billig, Y.:
{Irreducible representations for toroidal Lie algebras.}
J.Algebra {\bf 221}, 188-231 (1999).

\noindent
[B] Billig, Y., 
{A category of modules for the full toroidal Lie algebra.}  
Int.Math.Res. Not.  2006, Art. ID 68395, 46 pp.

\noindent
[BZ] Billig, Y., Zhao, K., 
{Weight modules over exp-polynomial Lie algebras.}
J.Pure Appl.Algebra  {\bf 191}  (2004),  no. 1-2, 23--42.

\noindent
[DLM] Dong, C., Li, H., Mason, G.,
{Vertex Lie algebras, vertex Poisson algebras and vertex algebras.}
in ``Recent developments in infinite-dimensional Lie algebras and conformal field theory'' 
(Charlottesville, VA, 2000),  69--96, Contemp. Math., {\bf 297}, Amer.Math.Soc., Providence, RI, 2002.

\noindent
[E1] Eswara Rao, S.,
{Irreducible representations of the Lie-algebra of the diffeomorphisms 
of a $d$-dimensional torus.}
J. Algebra {\bf 182} (1996), no. 2, 401--421. 

\noindent
[E2] Eswara Rao, S.,
{Partial classification of modules for Lie algebra of diffeomorphisms of $d$-dimensional torus.}
J. Math. Phys. {\bf 45} (2004), no. 8, 3322--3333. 

\noindent
[F] Futorny, V.M., 
{Irreducible non-dense $A_1^{(1)}$-modules.}
Pacific J. Math. {\bf 172} (1996), 83--99.

\noindent
[GKL] Gerasimov, A., Kharchev, S., Lebedev, D.,
{ Representation theory and quantum inverse scattering method: the open Toda chain and the hyperbolic Sutherland model.}
Int.Math.Res.Not.  2004,  no. 17, 823--854. 

\noindent
[K1] Kac, V.:
{\it Infinite dimensional Lie algebras.}
Cambridge: Cambridge University Press, 3rd edition, 1990.

\noindent
[K2]  Kac, V.:
{\it Vertex algebras for beginners.}
Second Edition, University Lecture Series, {\bf 10}, A.M.S., 1998.

\noindent
[KM] Khomenko, A., Mazorchuk, V., 
{On the determinant of Shapovalov form for generalized Verma modules. }
J.Algebra  {\bf 215}  (1999),  no. 1, 318--329.

\noindent
[LW] Lepowsky, J., Wilson, R.L.,
{Construction of affine Lie algebra $A_1^{(1)}$, }
Comm. Math. Phys. 62, 43-53 (1978).
 
\noindent
[L]  Li, H.,
{Local systems of vertex operators, vertex superalgebras and modules.}
{J.Pure Appl.Algebra} {\bf 109}, 143-195 (1996).

\noindent
[M]  Mathieu, O., 
{Classification of Harish-Chandra modules over the Virasoro algebra.} Invent. Math. {\bf 107} (1992), 225-234.

\end